\pdfoutput=1
\RequirePackage{ifpdf}
\ifpdf 
\documentclass[pdftex]{sigma}
\else
\documentclass{sigma}
\fi

\usepackage{dsfont}

\newcommand{\Span}{\operatorname{Span}}

\newcommand{\colu}{\underline{C}_p(\Gamma) }

\newcommand{\SL}{\operatorname{SL}}

\newcommand{\GL}{\operatorname{GL}}
\newcommand{\PGL}{\operatorname{PGL}}
\newcommand{\End}{\operatorname{End}}

\newcommand{\Mod}{\operatorname{Mod}}
\newcommand{\Add}{\operatorname{Add}}

\newcommand{\Spg}{\operatorname{Sp}( 2g, \mathbb{Z}/2\mathbb{Z})}

\newcommand{\HSg}{H_1(\Sigma_g, \mathbb{Z}/2\mathbb{Z})}
\newcommand{\cHHg}{\mathbb{C} [ H_1(H_g, \mathbb{Z}/2\mathbb{Z}) ] }
\newcommand{\cHSg}{\mathbb{C} [ H_1(\Sigma_g, \mathbb{Z}/2\mathbb{Z})] }
\newcommand{\mcg}{\Mod(\Sigma_g)}

\newcommand{\emcg}{\widetilde{\Mod}(\Sigma_g)}

\newcommand{\eemcg}{\widetilde{\Mod}(H_g)}

\newcommand{\quotient}[2]{{\raisebox{.2em}{$#1$}\left/\raisebox{-.2em}{$#2$}\right.}}

\newcommand{\eyesgraph}[3]{
\tikz[scale=0.3,baseline=-0.5ex]{
\draw (0,0) circle (1);
\draw (4,0) circle (1);
\draw (1,0) -- (3,0);
\draw (0,1) node[above]{#1};
\draw (2,0) node[above]{#2};
\draw (4,1) node[above]{#3};
}}

\newcommand{\Eyesgraph}{
\tikz[scale=0.3,baseline=-0.5ex]{
\draw (0,0) circle (1);
\draw (4,0) circle (1);
\draw (1,0) -- (3,0);
}
}

\newcommand{\thetagraph}[3]{
\tikz[scale=0.3,baseline=-0.5ex]{
\draw (0,0) circle (2);
\draw (0,2)--(0,-2);
\draw (-2,0) node[left]{#1};
\draw (0.2,0) node[left]{#2};
\draw (2,0) node[right]{#3};
}
}

\newcommand{\Tetahedrongraph}[6]{
\tikz[scale=0.4,baseline=-0.5ex]{
\draw (0,0) circle (2);
\draw (0,0)--(0,-2);
\draw (-1.35,1.35)--(0,0);
\draw (1.35,1.35)--(0,0);
\draw (0,2) node[above]{#1};
\draw (-1,1) node[below]{#2};
\draw (1,1) node[below]{#3};
\draw (-1.5,-1.3) node[left]{#4};
\draw (0,-1.3) node[left]{#5};
\draw (1.5,-1.3) node[right]{#6};
}
}

\newcommand{\invTetahedrongraph}[6]{
\tikz[scale=0.4,baseline=-0.5ex]{
\draw (0,0) circle (2);
\draw (0,0)--(0,2);
\draw (-1.35,-1.35)--(0,0);
\draw (1.35,-1.35)--(0,0);
\draw (0,-2) node[below]{#5};
\draw (-1,-1) node[above]{#4};
\draw (1,-1) node[above]{#6};
\draw (-1.5,1.3) node[left]{#1};
\draw (0,1.3) node[left]{#2};
\draw (1.5,1.3) node[right]{#3};
}
}

\newcommand{\Thetagraph}{
\tikz[scale=0.3,baseline=-0.5ex]{
\draw (0,0) circle (2);
\draw (0,2)--(0,-2);
}
}

\newcommand{\sixj}[6]{
\left\{
\begin{array}{ccc}
#1 & #2 & #3 \\ #4 & #5 & #6
\end{array}
\right\}
}

\newcommand{\muu}[1]{\mu_{#1}}

\newcommand{\eIgraph}[1]{
\tikz[scale=0.3,baseline=2ex]{
\draw (0,0) -- (1,1);
\draw (2,0) -- (1,1);
\draw (1,1) -- (1,2);
\draw (1,2) -- (0,3);
\draw (1,2) -- (2,3);
\draw (1,1.5) node[right]{#1};
\draw (0,0) node[left]{$a$};
\draw (0,3) node[left]{$b$};
\draw (2,3) node[right]{$c$};
\draw (2,0) node[right]{$d$};
}}

\newcommand{\eHgraph}[1]{
\tikz[scale=0.3,baseline=2ex]{
\draw (0,0) -- (1,1);
\draw (0,2) -- (1,1);
\draw (1,1) -- (2,1);
\draw (2,1) -- (3,0);
\draw (2,1) -- (3,2);
\draw (1.5,1) node[above]{#1};
\draw (0,0) node[left]{$a$};
\draw (0,2) node[left]{$b$};
\draw (3,2) node[right]{$c$};
\draw (3,0) node[right]{$d$};
}}

\usepackage{tikz, tikz-cd}
\numberwithin{equation}{section}

\newtheorem{Theorem}{Theorem}[section]
\newtheorem{Lemma}[Theorem]{Lemma}
\newtheorem{Proposition}[Theorem]{Proposition}
\newtheorem{Conjecture}[Theorem]{Conjecture}
 { \theoremstyle{definition}
\newtheorem{Remark}[Theorem]{Remark} }

\begin{document}

\allowdisplaybreaks

\newcommand{\arXivNumber}{1406.4389}

\renewcommand{\PaperNumber}{011}

\FirstPageHeading

\ShortArticleName{Decomposition of some Witten--Reshetikhin--Turaev Representations}

\ArticleName{Decomposition of some Witten--Reshetikhin--Turaev\\ Representations into Irreducible Factors}

\Author{Julien KORINMAN}

\AuthorNameForHeading{J.~Korinman}

\Address{Funda\c{c}\~ao Universidade Federal de S\~ao Carlos, Departamento de Matem\'atica,\\
Rod. Washington Lu\'{\i}s, Km 235, C.P.~676, 13565-905 S\~ao Carlos, SP, Brasil}
\Email{\href{mailto:julien.korinman@gmail.com}{julien.korinman@gmail.com}}
\URLaddress{\url{https://sites.google.com/site/homepagejulienkorinman/}}

\ArticleDates{Received October 29, 2017, in final form January 30, 2019; Published online February 12, 2019}

\Abstract{We decompose into irreducible factors the ${\rm SU}(2)$ Witten--Reshetikhin--Turaev representations of the mapping class group of a genus $2$ surface when the level is $p=4r$ and $p=2r^2$ with $r$ an odd prime and when $p=2r_1r_2$ with $r_1$, $r_2$ two distinct odd primes. Some partial generalizations in higher genus are also presented.}


\Classification{57R56; 57M60}

\section{Introduction}
Witten constructed in \cite{Wi2} a family of $(2+1)$-dimensional topological quantum field theories (TQFTs) using path integrals and the Chern--Simons action which gives a three-dimensional interpretation of the Jones polynomial. Each of these TQFTs induces a projective finite-dimensional representation of the mapping class group $\mcg$ of a genus $g$ closed oriented surface~$\Sigma_g$. Reshetikhin and Turaev made a rigorous construction of these TQFTs~\cite{RT} using representations of quantum groups. In this paper we will follow the skein theoretical construction of~\cite{BHMV2, Li2} to define these representations.

The Witten--Reshetikhin--Turaev projective representations lift to linear representations of a~central extension $\emcg$ of $\mcg$
\begin{gather*} \rho_{p,g}\colon \ \emcg \rightarrow \GL(V_p(\Sigma_g)).\end{gather*}
Here $p\geq 6$ is an even integer indexing the representations, called the \textit{level}, and $V_p(\Sigma_g)$ is a~finite-dimensional complex vector space. These representations are equipped with an invariant Hermitian non-degenerate form $\langle \, ,\, \rangle_{p,g}$.

The goal of this paper is to decompose some of these representations into irreducible factors. There are only few results in that direction. In \cite{BHMV2}, the authors construct an explicit proper invariant submodule of $V_p(\Sigma_g)$, when $4$ divides $p$. Roberts proved in \cite{Ro} that $\rho_{p,g}$ is irreducible if $\frac{p}{2}$ is an odd prime. An immediate extension of his proof shows that the representations $\rho_{18,g}$ are irreducible. In \cite{AF}, Andersen and Fjelstad proved that for $p=24,36,60$, the $\emcg$-module $V_p(\Sigma_g)$ contains at least three invariant submodules. The author gave in \cite{Koju1} an explicit decomposition into irreducible factors of the modules $V_p(\Sigma_1)$, arising in genus one, for arbitrary level $p\geq 3$. Note that one can extend the definition of the Witten--Reshetikhin--Turaev representations to mapping class groups of punctured surfaces, which are indexed by some coloring of the punctures in addition to the level. Koberda and Santharoubane showed in~\cite{KoberdaSantharoubane17} that any representation of a punctured surface with at least one puncture colored by $1$ is irreducible.

Given $r_1$, $r_2$ two distinct odd primes, there exists an unique even integer $x=x(r_1,r_2)\in \{ 1, \ldots, r_1r_2-2\}$ such that $x$ verifies either
\begin{gather*} \begin{cases}
x\equiv -2 &\pmod{r_1}, \\
x\equiv 0 & \pmod{r_2},
\end{cases} \qquad \text{or}\qquad \begin{cases}
x\equiv 0 &\pmod{r_1}, \\
x\equiv -2 & \pmod{r_2}.
\end{cases}\end{gather*}

The main result of this paper is the following:
\begin{Theorem}\label{main_th}\quad
\begin{enumerate}\itemsep=0pt
\item[$1.$] The modules $V_{18}(\Sigma_g)$ are simple for $g\geq 2$.
\item[$2.$] If $r$ is an odd prime, then $V_{4r}(\Sigma_2)$ is the sum of two simple submodules.
\item[$3.$] If $r$ is an odd prime, then $V_{2r^2}(\Sigma_2)$ is simple.
\item[$4.$] If $r_1$, $r_2$ are two distinct odd primes such that either $r_1,r_2>37$ or the element $x$ defined above satisfies $3x>2r_1r_2-4$. Then $V_{2r_1r_2}(\Sigma_2)$ is simple.
\end{enumerate}
\end{Theorem}

The main obstruction to extend the above theorem to higher genus is that we need to control which $6j$-symbols vanish. In the last section we will state a conjecture concerning the vanishing $6j$-symbols. We will then prove that this conjecture implies a generalization of Theorem~\ref{main_th} in higher genus. The author verified numerically the conjecture for small levels from which we deduce the:
\begin{Theorem}\label{theorem_high_genus} \quad
\begin{enumerate}\itemsep=0pt
\item[$1.$] Each module $V_{12}(\Sigma_3)$, $V_{20}(\Sigma_3)$, $V_{28}(\Sigma_3)$, $V_{44}(\Sigma_3)$, $V_{52}(\Sigma_3)$ is the direct sum of two simple submodules.
\item[$2.$] For any $g\geq 3$, the modules $V_{30}(\Sigma_g)$, $V_{66}(\Sigma_g)$ are simple.
\end{enumerate}
\end{Theorem}

The proof of Theorem~\ref{theorem_high_genus} relies on a numerical computation of the $6j$-symbols up to level~$66$.
\begin{Remark}In \cite{BHMV2} some representations $\rho_{p,g}$ are also defined when $p$ is odd. They verify $\rho_{2p,g}\cong \rho_{p,g}\otimes \rho'_{2,g}$. In particular if an odd level $r$ is such that $V_{2r}(\Sigma_g)$ is simple, then $V_{r}(\Sigma_g)$ is also simple. So our theorems extend to ${\rm SO}(3)$ cases as well.
\end{Remark}

\section[Skein construction of the Witten--Reshetikhin--Turaev representations]{Skein construction of the Witten--Reshetikhin--Turaev\\ representations}\label{section2}
Following \cite{BHMV2}, we will briefly define the Witten--Reshetikhin--Turaev representations and fix some notations.

\subsection[The spaces $V_{p}(\Sigma_g)$]{The spaces $\boldsymbol{V_{p}(\Sigma_g)}$}\label{section2.1}

\begin{figure}[!h]\centering
\includegraphics[width=4cm]{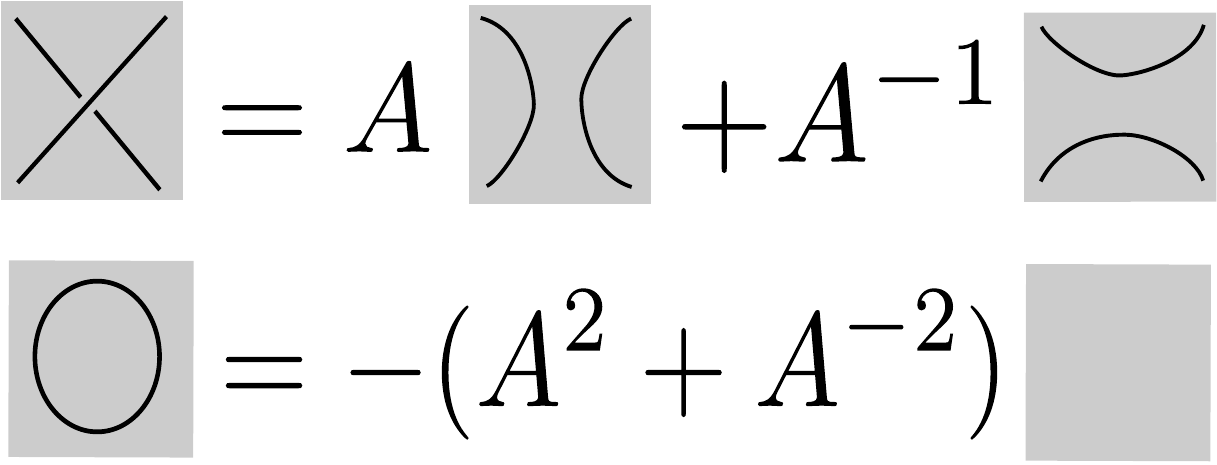}
\caption{The Kauffman-bracket skein relations defining the framed links invariant.}\label{fig_skein}
\end{figure}

Given an even integer $p\geq 6$, we denote by $A\in \mathbb{C}$ an arbitrary primitive $2p$-{th} root of unity. Given a compact oriented $3$-manifold $M$, a~\textit{framed link with $n$ components} $L\subset M$ is an isotopy class of an embedding of the disjoint union of $n$ copies of $S^1\times [0,1]$ into $M$. Using the Kauffman-bracket skein relation of Fig.~\ref{fig_skein}, we associate to any framed link $L\subset S^3$ an invariant $\langle L\rangle _p \in \mathbb{C}$.

Let $g\geq 1$ and denote by $\Sigma_g$ a compact oriented surface of genus $g$ and by $\Mod(\Sigma_g)$ the mapping class group of~$\Sigma_g$, namely the group of isotopy classes of orientation preserving homeomorphisms of $\Sigma_g$. Let $\mathcal{C}_g$ represents the set of isotopy classes of framed links (including the empty link) in an oriented genus $g$ handlebody~$H_g$. We fix a genus $g$ Heegaard splitting of the sphere, i.e., two homeomorphisms $S\colon \partial H_g \cong \overline{\partial H_g}$ and $ H_g\bigcup_{S\colon \partial H_g\rightarrow \overline{\partial H_g}} H_g \cong S^3$. For instance, define the handlebody $H_g=\big(D^2\times S^1\big)^{\# g}$ as a connected sums of $g$ copies of $D^2\times S^1$. Consider two oriented curves $L,M \subset \partial \big(D^2\times S^1\big)$ intersecting once positively. Denote by $L_i, M_i \subset \partial H_g$, $i=1, \ldots, g$ the images of $L$ and $M$ respectively in the $i$-th connected component of~$H_g$. One can choose the homeomorphism $S\colon \partial H_g \rightarrow H_g$ such that the image of the curve~$L_i$ is~$M_i$ for any $i\in \{1, \ldots, g\}$. The mapping class of $S$ is uniquely determined by this condition and we have a~homeomorphism $H_g\bigcup_S H_g \cong S^3$. We denote by $\varphi_1, \varphi_2 \colon H_g \hookrightarrow S^3$ the embeddings in the first and second factors.

Choose $L_1,L_2\in \mathcal{C}_g$. The above gluing defines a link $\varphi_1(L_1)\bigcup \varphi_2(L_2) \subset S^3$. The \textit{Hopf pairing} is the Hermitian form
\begin{gather*} ( \cdot , \cdot )_{g,p}^H \colon \ \mathbb{C}[\mathcal{C}_g]\times \mathbb{C}[\mathcal{C}_g] \rightarrow \mathbb{C}\end{gather*}
defined by \begin{gather*} ( L_1, L_2 t)_{g,p}^H := \left\langle \varphi_1(L_1)\bigcup \varphi_2(L_2) \right\rangle_p.\end{gather*}

Next we define the spaces $V_{p}(\Sigma_g)$ as the quotients:
\begin{gather*} V_{p}(\Sigma_g) := \quotient{ \mathbb{C}[\mathcal{C}_g]}{\ker \big(( \cdot , \cdot )_{g,p}^H \big)}.\end{gather*}

\looseness=-1 It is proved in \cite{BHMV2} that the vector spaces $V_{p}(\Sigma_g)$ are finite dimensional. Let us provide an explicit basis as follows. Given $g\geq 2$, choose a trivalent banded graph $\Gamma \subset H_g$ such that~$H_g$ retracts on $\Gamma$ by deformation. By banded graph we mean a thickening of the graph by an oriented surface. If $g=1$, $\Gamma $ represents the band $S^1\times \big[{-}\frac{1}{2}, \frac{1}{2}\big] \subset S^1\times D^2 \cong H_1$. We denote by $E(\Gamma)$ the set of edges of $\Gamma$. Set $I_p:=\big\{0, \ldots, \frac{p-4}{2}\big\}$, which will be called the set of \textit{colors} at level~$p$.

A triple of colors $i,j,k \in I_p$ is said \textit{$p$-admissible} if:
\begin{enumerate}\itemsep=0pt
\item[1)] $|i-j|\leq k \leq i+j$,
\item[2)] $i+j+k$ is even and $i+j+k\leq p-4$.
\end{enumerate}

A map $\sigma \colon E(\Gamma) \rightarrow I_p$ is a \textit{$p$-admissible coloring} of $\Gamma$ if for every three edges $e_1,e_2,e_3\in E(\Gamma)$ adjacent to a vertex, the triple $(\sigma(e_1),\sigma(e_2),\sigma(e_3))$ is $p$-admissible.

\begin{figure}[!h]\centering
\includegraphics[width=4.5cm]{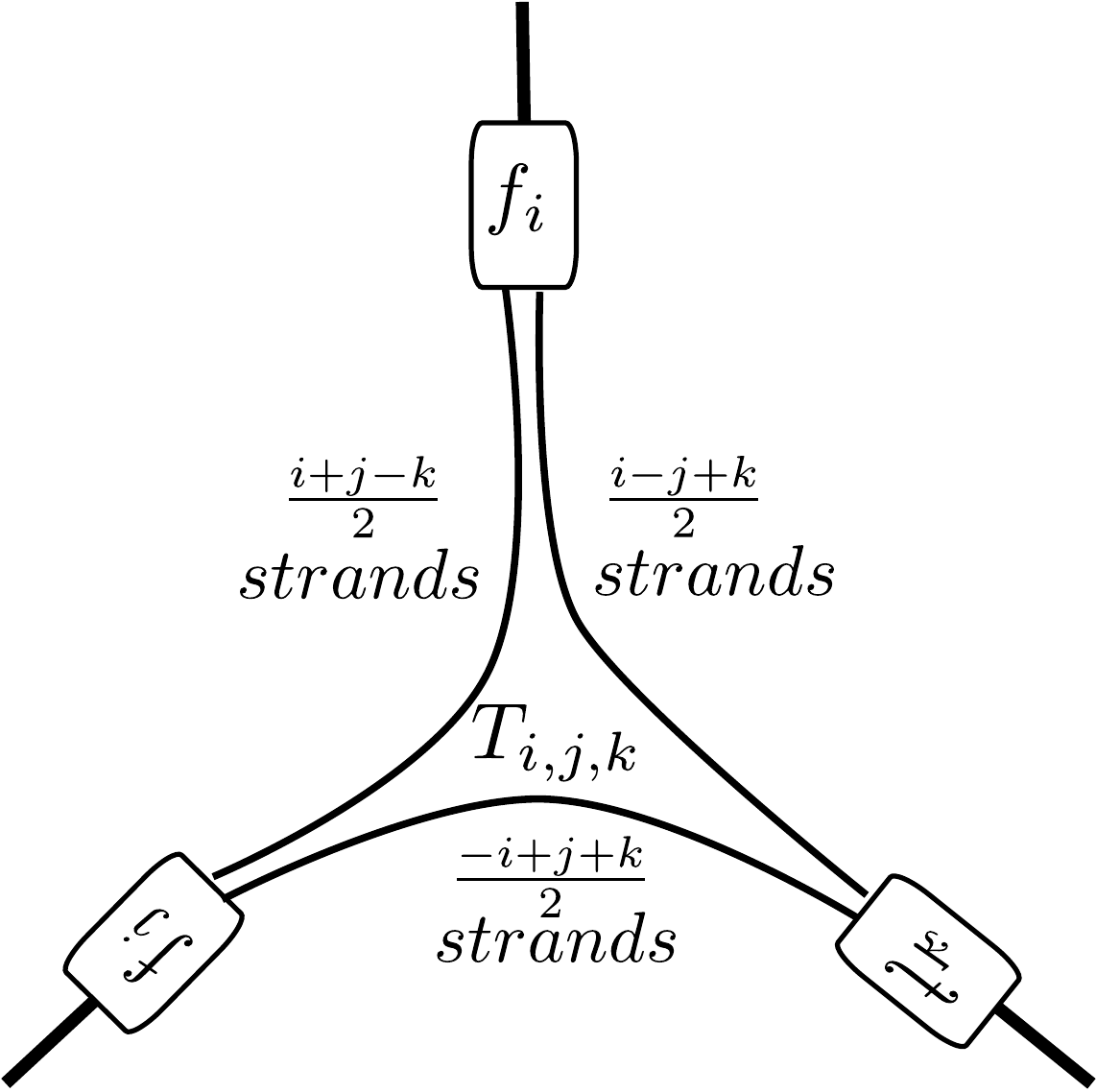}
\caption{The link $T_{i,j,k}$ used to connect three idempotents $f_i$, $f_j$ and $f_k$. The number above each three arcs denotes the number of parallel copies of the arc used to define the framed link.}\label{vertex_tangle}
\end{figure}

In \cite{Jones_proj, Wenzl}, the authors defined some idempotents $\big\{ f_0, \ldots, f_{\frac{p-4}{2}} \big\}$ of the Temperley--Lieb algebra with coefficient in $\mathbb{Q}(A)$ called \textit{the Jones--Wenzl idempotents}. To every $p$-admissible coloring $\sigma$ of~$\Gamma$, we associate a vector $u_{\sigma}\in V_{p}(\Sigma_g)$ as follows. Replace each edge $e\in E(\Gamma)$ by the Jones--Wenzl idempotent $f_{\sigma(e)}$. If $(e_1,e_2,e_3)$ are three edges adjacent to a vertex of~$\Gamma$, we connect the adjacent idempotents using the link $T_{\sigma(e_1),\sigma(e_2),\sigma(e_3)}$ drawn in Fig.~\ref{vertex_tangle}.

Theorem~4.11 of \cite{BHMV2} asserts that the vectors $u_{\sigma}$, where $\sigma$ belongs to the set of $p$-admissible colorings of $\Gamma$, form a basis of~$V_{p}(\Sigma_g)$.

\begin{figure}[!h]\centering
\centerline{\includegraphics[width=5cm]{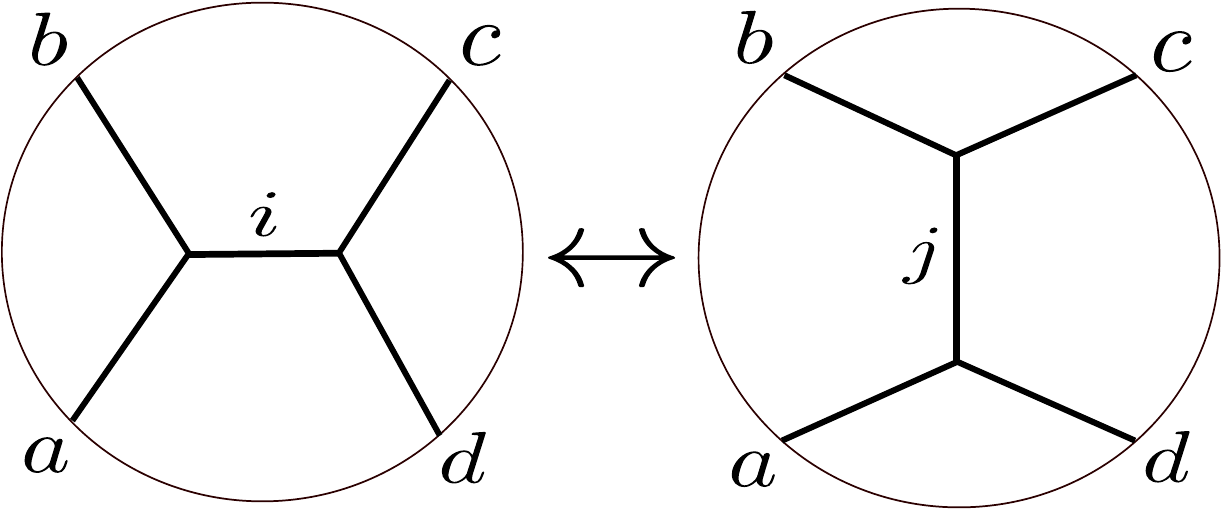}}
\caption{The two graphs $\Gamma_1$ on the left and $\Gamma_2$ on the right differ by a local Whitehead move.}\label{whitehead_move}
\end{figure}

The basis $u_{\sigma}$ depends on the choice of the embedded banded trivalent graph. One can transform any trivalent graph into any other one by a sequence of Whitehead moves and twists. We say that two banded trivalent graphs $\Gamma_1$ and $\Gamma_2$ embedded in a handlebody differ by a~Whitehead move if there exists a ball~$B^3$ intersecting each~$\Gamma_i$ transversally in four edges such that the two banded graphs coincide outside the ball and such that their intersection with $B^3$ is as drawn in Fig.~\ref{whitehead_move}.

Fix a $p$-admissible coloring of the graph outside $B^3$ with the four edges intersecting $B^3$ colored by $a$, $b$, $c$ and $d$. Denote by $\eHgraph{$i$}$ and $\eIgraph{$j$}$ the vectors associated to the coloring of~$\Gamma_1$ and~$\Gamma_2$ respectively with the edge inside $B^3$ colored by $i$ and $j$ respectively.

\begin{Lemma}[fusion-rule \cite{MV}]\label{fusion_rule}The vector $\eHgraph{i}$ belongs to the sub-space spanned by the vectors $\eIgraph{j}$ and decomposes as follows
\begin{gather*}\eHgraph{i} = \sum_j \sixj{a}{b}{j}{c}{d}{i} \eIgraph{j},\end{gather*}
where the sum runs through $p$-admissible colorings.
\end{Lemma}

In this formula, the coefficient $\sixj{a}{b}{j}{c}{d}{i}$ only depends on the colors $a$, $b$, $c$, $d$, $i$ and $j$ and is called \textit{recoupling coefficient} or \textit{$6j$-symbol}. We refer to \cite{MV} for a proof and an explicit computation of these coefficients. It follows from the formulas in \cite{MV} that if one of the colors $a$, $b$, $c$, $d$, $i$ or~$j$ is null, then the recoupling coefficient is nonzero.

\subsection{The Witten--Reshetikhin--Turaev representations}\label{section2.2}

We fix an orientation preserving homeomorphism $\alpha\colon \Sigma_g \rightarrow \partial H_g$.

Let $\phi \in \mcg$ be a mapping class associated to a homeomorphism of $\Sigma_g$ which extends to a homeomorphism $\widetilde{\phi}\colon H_g\rightarrow H_g$ through $\alpha$, that is such that $\phi=\alpha\circ \widetilde{\phi}_{| \partial H_g}\circ \alpha^{-1}$. Then $\widetilde{\phi}$ acts on $C_g$ and preserves the kernel of the Hopf pairing so acts on $V_{p}(\Sigma_g)$ by passing to the quotient. Denote by $\tilde{\rho}_{p,g}(\phi) \in \GL(V_{p}(\Sigma_g))$ the resulting operator. For instance, the Dehn twist along any curve $\gamma\subset \Sigma_g$ whose image $\alpha(\gamma)$ is contractible in~$H_g$, extends through~$\alpha$.

Let $\phi\in \mcg$ be the mapping class of a homeomorphism which extends to $H_g$ through $\alpha\circ S$. This extension also defines, by quotient, an operator on $V_{p}(\Sigma_g)$. We denote by $\tilde{\rho}_{p,g}(\phi)$ the dual of this operator for the Hopf pairing.

The elements of $\mcg$ which extend to $H_g$ either through $\alpha$ or through $\alpha\circ S$, generate the whole group $\mcg$. It is a non trivial fact that the associated operators $\tilde{\rho}_{p,g}(\phi)$ generate a projective representation:
\begin{gather*} \tilde{\rho}_{p,g} \colon \mcg \rightarrow \PGL(V_{p}(\Sigma_g)).\end{gather*}

We consider a central extension $\emcg$ of $\mcg$ that lifts the above projective representations to linear ones (see~\cite{GM13, MR}):
\begin{gather*} \rho_{p,g} \colon \ \emcg \rightarrow \GL(V_{p}(\Sigma_g)).\end{gather*}
These are the so-called Witten--Reshetikhin--Turaev representations.

The vector space $V_p(\Sigma_g)$ admits a non-degenerate Hermitian form, denoted $\langle \cdot, \cdot\rangle_{p,g}$ and called the \textit{invariant form}, which is invariant under the action of~$\emcg$. Moreover, for any trivalent banded graph, the basis associated to its $p$-admissible coloring is orthogonal for this form. The invariant form $\left<\cdot, \cdot\right>_{p,g}$ is related to the Hopf pairing $ (\cdot, \cdot )_{p,g}$ by the formula $ (v_1, v_2 )_{p,g}= \langle v_1, \rho_{p,g}(S)v_2\rangle _{p,g}$ for any $v_1,v_2\in V_p(\Sigma_g)$.

For each edge $e\in E(\Gamma)$ consider a disc $D_e$ properly embedded in $H_g$ that intersects $\Gamma$ transversely once in a point of the edge $e$. Note that the boundary curves $\gamma_e := \partial D_e \subset \partial H_g \cong \Sigma_g$ form a pants decomposition of $\Sigma_g$. Here and henceforth $T_e\in \mcg$ denotes the Dehn twist along $\gamma_e$.

From a classical property of the Jones--Wenzl idempotents, the authors of~\cite{BHMV2} proved that
\begin{gather*} \tilde{\rho}_{p,g}(T_e) \cdot u_{\sigma} = \muu{\sigma(e)} u_{\sigma}, \end{gather*}
where $\muu{i} := (-1)^i A^{i(i+2)}$ for every $i\in I_p$.

We fix the lift of $T_e$ in $\emcg$, still denoted $T_e$, so that $ \rho_{p,g}(T_e) \cdot u_{\sigma} = \muu{\sigma(e)} u_{\sigma}$.

We also fix the lift $S\in\emcg$ so that the matrix of $\rho_{p,g}(S)$ is the matrix of the Hopf pairing $(\cdot, \cdot )^H_{p,g}$ multiplied by the scalar $\eta:=\frac{1}{2p}(A\kappa)^3\big(A^2-A^{-2}\big)\sum\limits_{m=1}^{2p}{(-1)^mA^{-m^2}}\in \mathbb{C}^*$ where $\kappa$ satisfies $\kappa^6=A^{-6-\frac{p(p+1)}{2}}$. The scalar $\eta$ represents the sphere invariant in TQFT. We refer to~\cite{BHMV2} for a detailed discussion on $\eta$.

 Note that $S$ and the $ T_e, e\in E(\Gamma)$ generate $\emcg$ for any trivalent banded graph. Therefore, the above formulas determine uniquely the linear representation $\rho_{p,g}$.

Another description is derived as follows. Recall that $\mathcal{C}_g$ represents the set of isotopy classes of framed links in $H_g$ and denote by $\mathcal{K}_g$ the set of isotopy classes of framed links in $\Sigma_g\times [0,1]$. The vector space $\mathbb{C}[\mathcal{K}_g]$ has a natural algebra structure whose product is given by gluing two copies of $\Sigma_g\times[0,1]$. Moreover the homeomorphism $\alpha$ defines a structure of $\mathbb{C}[\mathcal{K}_g]$ left-module on $\mathbb{C}[\mathcal{C}_g]$ by gluing $\Sigma_g\times[0,1]$ to $H_g$. By passing to the quotient, we obtain a surjective map
\begin{gather*} \Add_p \colon \ \mathbb{C}[\mathcal{K}_g] \rightarrow \End(V_p(\Sigma_g)).\end{gather*}

If $K$ is a framed knot in an oriented $3$ manifold and $n\geq 0$, we note $K^n$ the framed link made of $n$ parallel copies of $K$. Recall that the framing is defined as a thickening of the link in an orientable surface. Pushing the link along the direction normal to the thickening defines the notion of parallelism. If $P(X)=\sum_i a_i X^i \in \mathbb{C}[X]$ is a polynomial and $K\in \mathcal{C}_g$ (resp $K\in \mathcal{K}_g)$ is a framed knot, we define $K(P):=\sum_i a_i K^i $ in $\mathbb{C}[\mathcal{C}_g]$ (resp in $\mathbb{C}[\mathcal{K}_g]$) and we call $K(P)$ the \textit{framed knot $K$ colored by $P$}.

Consider $\gamma\subset \Sigma_g$ a non-contractible oriented simple closed curve. Define $\gamma^+\subset \Sigma_g\times [0,1]$ the framed knot defined by $\gamma\subset \Sigma_g\times \big\{ \frac{1}{2} \big\}$ endowed with the framing given by a normal vector field making one turn in the counter-clockwise direction when circling once along $\gamma$ following the orientation. Define $\omega(X)=\sum\limits_{i=0}^ {\frac{p-4}{2}} (-1)^i [i+1]S_i(X) \in \mathbb{C}[X]$, where $S_i(X)$ represents the $i$-{th} Chebyshev polynomial of second kind defined by $S_0(X)=1$, $S_1(X)=X$ and $S_{i+2}(X)=XS_{i+1}(X)-S_i(X) $ and $[n]:= \frac{A^{2n}-A^{-2n}}{A^2 - A^{-2}}$. Then an important property of the Witten--Reshetikhin--Turaev representations is that the operator $\rho_{p,g}(T_{\gamma})$, associated to the lift of the Dehn twist along $\gamma$, is equal to the operator $\Add_p(\gamma^+(\omega))$ associated to the coloring by $\omega$ of $\gamma^+$. Since lifts of Dehn twists generate $\emcg$, this property gives an alternative definition of $\rho_{p,g}$.

\section{Cyclicity of the vacuum vector}

The vector $v_0\in V_p(\Sigma_g)$, which is the image of the empty link in $H_g$, will be called the \textit{vacuum vector} in genus~$g$. Denote by $\mathcal{A}_{p,g}$ the subalgebra of $\End(V_p(\Sigma_g))$ generated by the operators $\rho_{p,g}(\phi)$ for $\phi \in \emcg$. The key ingredient to prove Theorem~\ref{main_th} is to show that the vacuum vector is cyclic, i.e., that $\mathcal{A}_{p,g}\cdot v_0 =V_p(\Sigma_g)$.

\subsection{Decomposition in genus one of the Weil representations}
In \cite[Corollary 1.2]{Koju1}, the author gave an explicit decomposition of the genus one representations into irreducible factors which will be summarized next. We will derive from this decomposition the following:
\begin{Lemma}\label{lemma_genusone1} If $p=2r^2$ or $p=4r$ with $r$ an odd prime, or $p=2r_1r_2$ with $r_1$, $r_2$ distinct odd primes, then $v_0 \in V_p(\Sigma_1)$ is cyclic.
\end{Lemma}

The genus $1$ Weil representations at level $p$, are projective representations $\pi_p \colon \SL_2(\mathbb{Z}) \rightarrow \PGL(U_p)$ where $U_p$ is a $p$-dimensional complex vector space with a canonical basis $\{e_i, \, i \in \mathbb{Z}/p\mathbb{Z} \}$. They are defined on the generators $T:= \left(\begin{smallmatrix} 1&1\\0&1 \end{smallmatrix}\right)$ and $S:=\left(\begin{smallmatrix} 0& -1 \\ 1&0 \end{smallmatrix}\right)$ by the projective classes of the operators:
 \begin{gather*}
\pi_p (S) = \frac{1}{\sqrt{p}} \big( A^{-ij}\big)_{i,j \in \mathbb{Z}/p\mathbb{Z}},\qquad
\pi_p(T) = \big( A^{i^2}\delta_{i,j} \big)_{i,j \in \mathbb{Z}/p\mathbb{Z}}.
\end{gather*}

 Here the level is an integer $p\geq 2$ not necessarily even. When $p$ is even, $A$ is a primitive $2p$-{th} root of unity. When $p$ is odd, $A$ represents a primitive $p$-{th} root of unity. The decomposition into irreducible factors of the genus one Weil representations is the following:

\begin{Theorem}[{\cite[Theorem ~1.1]{Koju1}}] \label{th_weil}\quad
\begin{enumerate}\itemsep=0pt
\item[$1.$] If $a$ and $b$ are coprime, then $U_{ab}\cong U_a \otimes U_b$.
\item[$2.$] If $r$ is prime and $n\geq 1$, then there exists some module $W_{r^{n+2}}$ such that $U_{r^{n+2}}\cong U_{r^n}\oplus W_{r^{n+2}}$.
\item[$3.$] If $r$ is an odd prime, then $U_{r^2}\cong \mathds{1}\oplus W_{r^2}$ where $\mathds{1}$ is the one-dimensional trivial representation.
\item[$4.$] Each one of the modules $U_p$ and $W_{r^n}$ is the direct sum of two submodules $U_p\cong U_p^-\oplus U_p^+$, $W_{r^n}\cong W_{r^n}^+\oplus W_{r^n}^-$.
\item[$5.$] Every module of the form $B_{1}\otimes \dots \otimes B_{k}$, where $B_{i}$ is either $U_{r}^+$, $U_{r}^-$, $U_2$, $U_4^+$, $U_4^-$, $W_{r^n}^+$ or $W_{r^n}^-$ for $r$ prime and where the $B_i$ have pairwise coprime levels, is simple.
\end{enumerate}
\end{Theorem}
Given $n\geq0$ and $N\geq 1$ two non-negative integers, denote by $[n]_N$ the class of $n$ in $\mathbb{Z}/N\mathbb{Z}$. The isomorphism $U_{ab}\rightarrow U_a\otimes U_b$ in Theorem~\ref{th_weil}(1) sends the vector $e_{[n]_{ab}}$ to $e_{[n]_a}\otimes e_{[n]_b}$. The submodules $U_p^{\pm}\subset U_p$ arising in Theorem~\ref{th_weil}(3) are spanned by the vectors $e_i^{\pm}:= e_i \pm e_{-i}$.

The decomposition in simple submodules of $V_p(\Sigma_1)$ then follows from the above theorem and from the fact that if $p=2r\geq 6$, the map $\Psi\colon U_p^- \rightarrow V_{p}(\Sigma_1)$
defined by $\Psi(e_i^-) = u_{i+r-1}$ is an isomorphism of $\SL_2(\mathbb{Z})$-projective modules. This was proved in \cite{FK} (see also \cite[Theorem~5.1]{Koju1}).

\begin{proof}[Proof of Lemma \ref{lemma_genusone1}]If $p=2r^2$ with $r$ an odd prime, Theorem~\ref{th_weil}(5) implies that $U_{2 r^2}^- \cong V_{2r^2}(\Sigma_1)$ is simple, hence $v_0$ is cyclic. Assume next that $p=2r_1r_2$ with $r_1$, $r_2$ two distinct odd primes. Theorem~\ref{th_weil} provides an isomorphism of $\SL_2(\mathbb{Z})$ projective modules
\begin{gather*} \theta \colon \ U_{2r_1r_2}^- \cong \big(U_2\otimes U_{r_1}^-\otimes U_{r_2}^+\big) \oplus \big(U_2\otimes U_{r_1}^+ \otimes U_{r_2}^-\big).\end{gather*}
Observe that both modules $U_2\otimes U_{r_1}^-\otimes U_{r_2}^+$ and $U_2\otimes U_{r_1}^+ \otimes U_{r_2}^-$ are simple. To prove that $v_0 \in V_{2r_1r_2}(\Sigma_1)$ is cyclic, we need to show that $v:= \theta\circ (\Psi)^{-1} (v_0) $ has non trivial projection in both submodules. First we have $\Psi^{-1}(v_0) = e_{-1} - e_1$ (recall that indices are considered modulo~$p$). Then we compute:
\begin{gather*} v = e_0\otimes e_{-1} \otimes e_{-1} - e_0\otimes e_1\otimes e_1 = \left( \tfrac{1}{2} e_0 \otimes e_{-1}^- \otimes e_{-1}^+ \right) + \left(\tfrac{1}{2} e_0\otimes e_{-1}^+\otimes e_{-1}^- \right).\end{gather*}
The above decomposition shows that projections onto both submodules are non trivial, so~$v_0$ is cyclic.

The case $p=4r$ is similar. We start with the decomposition
\begin{gather*} \theta \colon \ U_{4r}^- \cong \big(U_4^- \otimes U_r^+\big) \oplus \big(U_4^+ \otimes U_r^-\big).\end{gather*}
Note that the modules on the right-hand side are simple. We keep the notation $v:=\theta\circ (\Psi)^{-1}(v_0)$ and, setting $n=2r-1$, we compute
\begin{gather*} v= e_{[n]_4}\otimes e_{[n]_r} - e_{[-n]_4}\otimes e_{[-n]_r} = e_1\otimes e_{-1} - e_{-1}\otimes e_1 \\
\hphantom{v}{} = \tfrac{1}{2}(e_1-e_{-1})\otimes (e_1+e_{-1}) - \tfrac{1}{2}(e_1+e_{-1})\otimes (e_1-e_{-1}).
\end{gather*}
So the projections of $v_0$ on the simple submodules of $U_p^-$ are nonzero. This concludes the proof.
\end{proof}

\subsection{From genus one to higher genus}

The fact that the vacuum vector is cyclic in genus one will give us information on the cyclic space of the vacuum vector in higher genus. We now describe the lemma that states this relation.

Choose $\Gamma_g \subset {S}^3$ a banded trivalent graph whose underlying graph has genus $g\geq 2$ embedded in the three-sphere. Consider $\Gamma_{g'}\subset \Gamma_g$ a sub-banded graph whose underlying graph has genus $g'\leq g$. Let~$H_g$ be a tubular neighborhood of $\Gamma_g$ and $H_{g'}^1$, $H_{g'}^2$ two tubular neighborhoods of~$\Gamma_{g'}$ such that $H_{g'}^1\subset H_g \subset H_{g'}^2$. Denote by $i_1\colon H_{g'}^1 \hookrightarrow H_g$ and $i_2\colon H_g \hookrightarrow H_{g'}^2$ the embeddings. The embedding $i_1$ induces a linear map $i_1^*\colon \mathbb{C}[\mathcal{C}_{g'}] \rightarrow \mathbb{C}[\mathcal{C}_g]$ sending a framed link $L\subset H_{g'}^1$ to $i_1(L) \subset H_g$. The morphism $i_1^*$ sends the kernel of the Hopf pairing in genus~$g'$ to a sub-space of the kernel of the Hopf pairing in genus~$g$. Hence it induces a linear map $j_1 \colon V_p(\Sigma_{g'}) \rightarrow V_p(\Sigma_g)$. Similarly, the embedding $i_2$ induces a linear map $j_2\colon V_p(\Sigma_g) \rightarrow V_p(\Sigma_{g'})$. Since $i_2 \circ i_1$ is a~retraction by deformation, the composition $j_2\circ j_1$ is the identity map of $V_p(\Sigma_{g'})$.

A framed link $L\subset H_g$ is called \textit{aligned} if there exist some oriented simple closed curves $\gamma_1, \ldots, \gamma_n\subset \Sigma_g$ such that $\gamma_1^+ \cdot \gamma_2^+ \cdots \gamma_n^+ \cdot \varnothing = L$, where we used the algebra structure of $\mathbb{C}[\mathcal{K}_g]$ and the left-module structure of $\mathbb{C}[\mathcal{C}_g]$ defined at the end of the previous section. Recall that the image of a Dehn twist $\rho_p(T_{\gamma})$ is equal to the operator $\Add_p(\gamma^+(\omega))$. It follows that the cyclic space of the vacuum vector is spanned by the vectors $[L(\omega)]$ obtained from an aligned framed link by coloring each of its connected components by $\omega$. By construction, the map $i_1^*$ sends aligned links to aligned links.

Fix $\sigma'$ a $p$-admissible coloring of $\Gamma_{g'}$. Consider the $p$-admissible coloring $\sigma$ of $\Gamma_g$ defined by
\begin{gather*} \sigma(e)= \begin{cases}
\sigma'(e), & \text{if $e$ is an edge of $\Gamma_{g'}\subset \Gamma_g$},\\
0, & \text{otherwise}.
\end{cases}\end{gather*}

\begin{Lemma}\label{lemma_reduction}If $u_{\sigma'}\in V_p(\Sigma_{g'})$ is in the cyclic space of the vacuum vector in genus~$g'$, then $u_{\sigma}\in V_p(\Sigma_g)$ also belongs to the cyclic space of the vacuum vector in genus~$g$.
\end{Lemma}

\begin{proof}[Proof of Lemma \ref{lemma_reduction}]
By hypothesis, there exists a linear combination of aligned framed links $\mathcal{L}\in \mathbb{C}[\mathcal{C}_{g'}]$ such that $[\mathcal{L}(\omega)]= u_{\sigma'} \in V_p(\Sigma_{g'})$. We now show that $[i_1^*(\mathcal{L})(\omega)]=u_{\sigma}\in V_p(\Sigma_g)$. Since~$i_1^*$ sends aligned framed links to aligned framed links, the claim will follow. Denote by \mbox{$T_{\sigma'}\in \mathbb{C}[\mathcal{C}_{g'}]$} the element obtained by replacing each edge $e$ of $\Gamma_g$ by the Jones--Wenzl idempotent~$f_{\sigma'(e)}$ and connecting the resulting elements as described in Section~\ref{section2.1}. By definition, we have $u_{\sigma'}=[T_{\sigma'}]$ and $u_{\sigma}=[j_1 ( T_{\sigma'})]$. To prove that the vectors $i_1^*(\mathcal{L})(\omega)$ and $j_1 ( T_{\sigma'})$ represent the same class in the quotient $V_p(\Sigma_g)$, we need to prove that their Hopf pairing with any framed link are equal. Let $K\subset S^3\setminus \mathring{H}_g$ be a framed link. We compute
\begin{gather*}
\big( i_1^*(\mathcal{L})(\omega), K \big)^H_{g, p} = ( \mathcal{L}(\omega), j_2(K) )^H_{g', p} = ( T_{\sigma'}, j_2(K) )^H_{g',p}= ( j_1(T_{\sigma'}), K )^H_{g,p},
\end{gather*}
where we passed from the first line to the second line by using the fact that $[\mathcal{L}(\omega)]=u_{\sigma'}= [T_{\sigma'}]\in V_p(\Sigma_{g'})$. This proves that $u_{\sigma}=[i_1^*(\mathcal{L})(\omega)]$ and concludes the proof.
\end{proof}

\subsection{Cyclicity in genus 2}
The goal of this subsection is to prove the following:
\begin{Proposition}\label{prop_cyclicity_high_genus}
Suppose that either $p=2r^2$ or $p=4r$ with $r$ an odd prime, or $p=2r_1r_2$ with $r_1$, $r_2$ two distinct odd primes. Then the vacuum vector $v_0\in V_p(\Sigma_2)$ is cyclic.
\end{Proposition}

Let $\Gamma\subset H_g$ a trivalent graph, as in Section~\ref{section2}. Two $p$-admissible colorings $\sigma_1$, $\sigma_2$ of $\Gamma$ are equivalent if $\mu_{\sigma_1(e)}=\mu_{\sigma_2(e)}$, for every edges $e \in E(\Gamma)$.

We denote by $\colu$ the set of equivalence classes of $p$-admissible colorings of~$\Gamma$. Given $[\sigma]\in \colu$, we associate the subspace
\begin{gather*} W_{[\sigma]} := \Span \{ u_{\sigma'},\, \sigma'\in [\sigma] \} \subset V_p(\Sigma_g). \end{gather*}

\begin{Lemma}\label{lemma1}If $X\subset V_p(\Sigma_g)$ is an invariant $\emcg$-submodule, then $ X=\bigoplus_{[\sigma]\in \colu} X\cap W_{[\sigma]}$.
\end{Lemma}

\begin{proof} The matrices $\rho_{p,g}(T_e)$, for $e\in E(\Gamma)$, generate a commutative subalgebra $\mathcal{T}$ of $\mathcal{A}_{p,g}$. To every character $\chi \colon \mathcal{T}\rightarrow \mathbb{C}^*$, we associate a subspace $V_p(\Sigma_g, \chi)$ formed by the vectors $v\in V_p(\Sigma_g)$ such that $\rho_{p,g}(T_e)v=\chi(\rho_{p,g}(T_e))v$. The set $\colu$ is in natural bijection with the characters $\chi$ such that $V_p(\Sigma_g,\chi)\neq 0$ and the spaces $W_{[\sigma]}$ correspond to the associated subspaces $V_p(\Sigma_g, \chi)$. Since the orthogonal projector on $X$ commutes with the elements of the algebra~$\mathcal{T}$, it preserves the subspaces~$W_{[\sigma]}$.
\end{proof}

The strategy to prove Proposition~\ref{prop_cyclicity_high_genus} is to apply Lemma~\ref{lemma1} to the subspace $ X:= ( \mathcal{A}_{p,2}\cdot v_0 ) ^{\bot}$ which is the orthogonal for the invariant form $\langle \cdot, \cdot \rangle_{p,g}$ of the cyclic space generated by the vacuum vector. Since the invariant form and the Hopf pairing are related by the formula $(v_1, v_2)_{p,g}=\langle v_1, \rho_{p,g}(S)\cdot v_2\rangle_{p,g}$, the space~$X$ is also the orthogonal of the cyclic space $\mathcal{A}_{p,2}\cdot v_0$ for the Hopf pairing.

\begin{figure}[!h]\centering
\includegraphics[width=4.8cm]{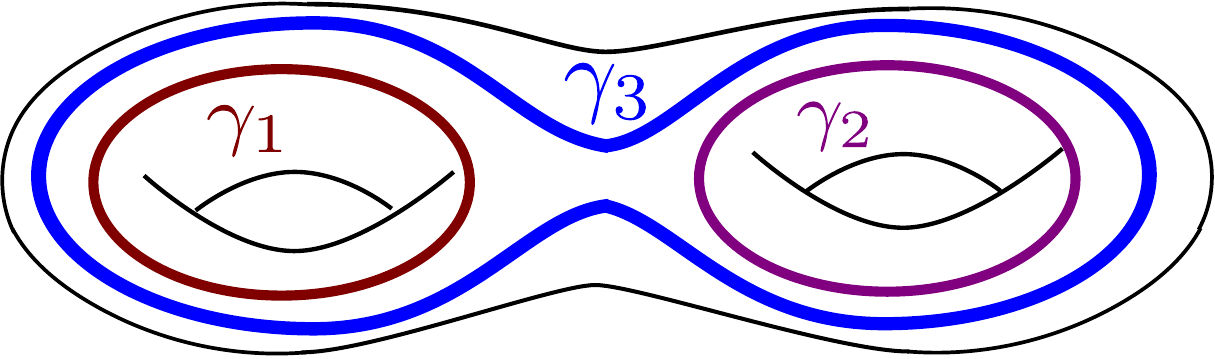}
\caption{The three framed knots defining the vectors $w_{a,b,c}$. The framing is defined by thickening the knots in a surface parallel to the boundary.}\label{troiscourbes}
\end{figure}

Let $\gamma_1, \gamma_2, \gamma_3 \subset H_2$ be the framed knots of Fig.~\ref{troiscourbes}. If $a$, $b$, $c$ are non-negative integers, define $w_{a,b,c}\in V_p(\Sigma_2)$ to be the class of the framed link made of~$a$ parallel copies of $\gamma_1$, $b$ parallel copies of $\gamma_2$ and~$c$ parallel copies of~$\gamma_3$.

\begin{Lemma}\label{lemma2}If $p=4r$, with $r$ an odd prime, or if $p=2r_1r_2$, with $r_1$, $r_2$ two distinct odd primes, then $ w_{a,b,c} \in \mathcal{A}_{p,2}\cdot v_0$ for any $a,b,c\in \{ 0, 1\}$.
\end{Lemma}

\begin{proof}Lemmas \ref{lemma_genusone1} and \ref{lemma_reduction} imply that $w_{1,0,0}$, $ w_{0,1,0}$ and $w_{1,1,0}$ belong to the cyclic space of the vacuum vector. It remains to show that $w_{1,1,1}$ also belongs to this space. Recall that the coefficient $x_k:=\sixj{2}{2}{k}{2}{2}{0}$ is nonzero when $(2,2,k)$ is $p$-admissible.
 Using Lemma~\ref{fusion_rule}, we have the following system
\begin{gather*}
\thetagraph{2}{0}{2} - x_0\eyesgraph{2}{0}{2} = x_2 \eyesgraph{2}{2}{2} +x_4 \eyesgraph{2}{4}{2}\,, \\
\rho_{p,2}(T_e) \cdot \thetagraph{2}{0}{2} - x_0\eyesgraph{2}{0}{2} = \muu{2} x_2 \eyesgraph{2}{2}{2} +\muu{4}x_4 \eyesgraph{2}{4}{2}\,,
\end{gather*}
where $T_e$ is (a lift of) the Dehn twist around the middle edge of the graph $\Eyesgraph$. Both vectors on the left-hand side belong to the cyclic space of $v_0$ by Lemmas \ref{lemma_genusone1} and \ref{lemma_reduction}.

Since the coefficients $x_2$ and $x_4$ are nonzero and $\mu_2\neq \mu_4$, the matrix $\left(\begin{smallmatrix} x_2 & x_4 \\ \mu_2 x_2 & \mu_4 x_4 \end{smallmatrix}\right)$ is invertible. We deduce from the above linear system that both vectors $\eyesgraph{2}{2}{2}$ and $\eyesgraph{2}{4}{2}$ belong the cyclic space of the vacuum vector. Since the vector~$w_{1,1,1}$ belongs to the space spanned by the vectors $\eyesgraph{2}{$k$}{2}$ for $k=0,2,4$, it belongs the cyclic space of~$v_0$.
\end{proof}

Recall the notation $\mu_i=(-1)^iA^{i(i+2)}$ and the equivalence relation $i\sim j$ if $\mu_i=\mu_j$ used to define the spaces~$W_{[\sigma]}$. The following lemma describes this equivalence relation.
\begin{Lemma}\label{lemmachoucroute} Let $i,j\in I_p$. Then:
\begin{enumerate}\itemsep=0pt
\item[$1.$] When $p=2r^2$ with $r$ an odd prime, then $\mu_i=\mu_j$ if and only if $i\equiv j\equiv -1\pmod{r}$ and~$i$,~$j$ have the same parity.
\item[$2.$] When $p=4r$ with $r$ an odd prime, then $\mu_i=\mu_j$ if and only if either $i=j$ or $(i=\frac{p-4}{2} - j$ and $i$ is even$)$.
\item[$3.$] When $p=2r_1r_2$ with $r_1$, $r_2$ distinct odd primes, then $\mu_i=\mu_j$ if and only if either $i=j$ or~$j$ is the only element satisfying either
\begin{gather*} \begin{cases}
i\equiv j & \pmod{2r_1}, \\
i\equiv -j-2 & \pmod{r_2}
\end{cases}
\qquad \text{or}\qquad
 \begin{cases}
i\equiv j & \pmod{2r_2}, \\
i\equiv -j-2 & \pmod{r_1}.
\end{cases}
\end{gather*}
\item[$4.$] When $p=18$, all $\mu_i$ are pairwise distinct.
 \end{enumerate}
\end{Lemma}

\begin{proof}The case $p=18$ is proved by a straightforward computation. For the other cases, first note that:
\begin{gather}\label{mui}
 \mu_i=\mu_j \Leftrightarrow A^{p(i+j)+(i-j)(i+j+2)}=1 \Leftrightarrow p(i+j)+(i-j)(i+j+2)\equiv 0 \pmod{2p}.
\end{gather}
When restricted modulo $4$, equation~\eqref{mui} implies that $i$ and $j$ have same parity.

When $p=4r$, equation \eqref{mui} implies $(i-j)(i+j+2)\equiv 0 \pmod{r}$. Since $i\neq j$ and $i$, $j$ have same parity, then $i=2r-2-j$. When restricted modulo $8$, equation~\eqref{mui} implies that $i\equiv j\pmod{4}$ or $i\equiv j \equiv 0 \pmod{2}$. The relations $i=2r-2-j$ and $i\equiv j \pmod{4}$ forbid $i$ and $j$ to be odd. Hence $i$ and $j$ are even and $i=2r-2-j$.

When $p=2r_1r_2$, equation~\eqref{mui} implies $(i-j)(i+j+2) \equiv 0\pmod{r_1r_2}$. Since $i\neq j$ and $i$, $j$ have same parity, we must have either $i\equiv j \pmod{r_1}$ and $i\equiv -j-2 \pmod{r_2}$ or $i\equiv j \pmod{r_2}$ and $i\equiv -j-2 \pmod{r_1}$.

Finally, when $p=2r^2$, equation~\eqref{mui} implies $(i-j)(i+j+2) \equiv 0\pmod{r^2}$. Since $i\neq j$ and $i$, $j$ have same parity, $r$ divides both $i-j$ and $i+j+2$ and thus $i\equiv j \equiv -1 \pmod{r}$.
\end{proof}

\begin{Lemma}\label{lemma3}If $p=2r^2$, with $r$ an odd prime, and $\sigma$ is a $p$-admissible coloring of $\Gamma=\Thetagraph$ such that $\sigma(e) \not\equiv -1 \pmod{r}$, for all $e\in E(\Gamma)$, then $u_{\sigma}\in \mathcal{A}_{p,2}\cdot v_0$. Moreover, if $0\leq a,b,c \leq \frac{r-3}{2}$, then $w_{a,b,c}\in \mathcal{A}_{p,2}\cdot v_0$.
\end{Lemma}

\begin{proof} When $\sigma$ satisfies the hypothesis of the lemma, Lemma~\ref{lemmachoucroute} implies that the sub\-spa\-ce~$W_{[\sigma]}$ is one dimensional. Lemma~\ref{lemma1} implies that this subspace is included either in $\mathcal{A}_{p,2}\cdot v_0$, or in its orthogonal. Moreover the Hopf pairing $\left( u_{\sigma}, v_0 \right)^H_{p,2}$ is nonzero since it is equal to a $p$-admissible $3j$-symbol. Hence we have proved the first statement of the lemma. In particular, we have the inclusion:
\begin{gather*} S:=\Span \left(\thetagraph{$u$}{$v$}{$w$}, 0\leq u,v,w \leq r-2 \right) \subset \mathcal{A}_{p,2}\cdot v_0.\end{gather*}

The proof of the second statement follows from the fact that the vector $w_{a,b,c}$ belongs to the subspace $S$ whenever we have $a+c\leq r-2$, $b+c\leq r-2$ and $a+b\leq r-2$.
\end{proof}

\begin{proof}[Proof of Proposition~\ref{prop_cyclicity_high_genus}]First assume that $p\neq 18$. Consider the graph $\Gamma=\Thetagraph$, a class $[\sigma]\in \colu$, and choose a vector
\begin{gather*}
 v \in W_{[\sigma]}\bigcap \left( A_{p,2}\cdot v_0\right)^{\bot},\qquad v = \sum_{\tau\in [\sigma]}{\alpha_{\tau}u_{\tau}}.
\end{gather*}

From Lemma~\ref{lemma1}, we must show that $v=0$ to conclude. To this purpose, we will find $\dim \left(W_{[\sigma]}\right)$ linearly independent equations satisfied by the coefficients $\alpha_{\tau}$.

Denote by $F\subset \mathbb{N}^{E(\Gamma)}$ the set of functions $f$ such that:
\begin{itemize}\itemsep=0pt
\item $f(e) \in \{0,1\}$, for all $e\in E(\Gamma)$, if $p=4r$ or $p=2r_1r_2$,
\item $f(e) \in \big\{ 0, \ldots, \frac{r-3}{2} \big\}$, for all $e \in E(\Gamma)$, if $p=2r^2$.
\end{itemize}
Write $E(\Gamma)=\{e_1, e_2, e_3\}$ and, given $f\in F$, we define the vector $w_f:= w_{f(e_1), f(e_2), f(e_3)}$. Lemmas~\ref{lemma2} and~\ref{lemma3} imply that $w_f\in \mathcal{A}_{p,2}\cdot v_0$ for all $f\in F$. By definition of $v$, we have that $(w_f, v)_{p,2}^H = 0$ for all $f\in F$. This implies that
\begin{gather*}
 \sum_{\tau\in [\sigma]}{\left( \prod_{e\in E(\Gamma)}{\lambda_{\tau(e)}^{f(e)}}\right) \alpha_{\tau} \left( u_{\tau}, v_0 \right)_{p,2}^H } =0, \qquad \text{for all} \quad f\in F,
\end{gather*}
where $\lambda_i = -\big(A^{2(i+1)}+A^{-2(i+1)}\big)$. Since the complex numbers $( u_{\tau}, v_0)_{p,2}^H$ are $p$-admissible~$3j$ symbols, they are nonzero. So it is enough to show that the matrix
\begin{gather*} M := \left(\prod_{e\in E(\Gamma)}{\lambda_{\tau(e)}^{f(e)}}\right)_{\substack{\tau\in[\sigma] \\ f\in F}} \end{gather*}
has independent lines.

We now define an invertible square matrix $\widetilde{M}$ such that $M$ is obtained from $\widetilde{M}$ by removing some lines. Given $i\in I_p$, we define the set
\begin{gather*} \omega(i) := \big\{ j\in I_p, \mbox{ so that }\muu{i}=\muu{j} \big\}.\end{gather*}
Denote by $\# \omega(i)$ its cardinal and define the Vandermonde matrix
\begin{gather*} N[i] := \big( \lambda_j^n \big) _{\substack{ j \in \omega(i)\\ 0\leq n \leq \# \omega(i) -1 }}.\end{gather*}
Since $\lambda_i\neq \lambda_j$ when $i\neq j $, the matrix $N[i]$ is invertible. The matrix $ \widetilde{M} := N[e_1]\otimes N[e_2] \otimes N[e_3]$ is invertible and $M$ is obtained from $\widetilde{M}$ by removing the lines corresponding to non $p$-admissible colorings of $\Gamma$. This ends the proof when $p\neq 18$.

The proof of the lemma when $p=18$ is similar to Roberts' proof in~\cite{Ro} which only relies on the fact that the coefficients $\mu_i$ are pairwise distinct. We briefly reproduce it. Let $K\subset H_1$ the framed knot $\{0\}\times S^1 \subset D^2\times S^1=H_1$ with trivial framing, so that $[K]=u_1\in V_p(\Sigma_1)$. The Hopf pairing of $K^i(\omega)$ with $u_j$ is $(\mu_j)^i$. Since the $\mu_j$ are pairwise distinct, the Vandermonde matrix $\big(\big( \big[K^i(\omega)\big], u_j \big)_{18, 1}^H\big)_{i,j}$ is invertible. Since the Hopf pairing is non degenerate, it follows that the vectors $[K^i (\omega)]$ for $i\in \{0, \ldots, 8 \}$ form a basis of $V_{18}(\Sigma_1)$. In particular, $u_1$ is a linear combination of vectors $[K^i(\omega)]$. Now choose $L\subset H_g$ an aligned link. Replacing each connected component of $L$ by the above linear combination of parallel copies colored by $\omega$, we see that the class $[L]\in V_{18}(\Sigma_g)$ is equal to the class of a linear combination of aligned links colored by~$\omega$, thus belongs to $\mathcal{A}_{18,g}\cdot v_0$. Since the vectors $[L]\in V_{18}(\Sigma_g)$, with $L$ aligned, span the whole space~$V_{18}(\Sigma_g)$, this concludes the proof.
\end{proof}

\section{Decomposition into irreducible factors}

 In this section, we will prove Theorem~\ref{main_th}. Denote by $(\mathcal{A}_{p,g})'$ the commutant of the algebra~$\mathcal{A}_{p,g}$, i.e., the subspace of $\End(V_p(\Sigma_g))$ of operators commuting with all operators $\rho_{p,g}(\phi)$ for $\phi\in \emcg$.

The dimension of $(\mathcal{A}_{p,g})'$ is equal to the number of simple submodules of $V_p(\Sigma_g)$. Thus we have to show that $\dim( (\mathcal{A}_{p,2})')$ is one if $p=2r^2$ and $p=2r_1r_2$ and two when $p=4r$. We will also prove that $\dim( (\mathcal{A}_{18,g})')=1$.

Consider the linear map $f\colon (\mathcal{A}_{p,g})' \hookrightarrow V_p(\Sigma_g)$ defined by $f(\theta)= \theta\cdot v_0$. Since $v_0$ is cyclic by Proposition~\ref{prop_cyclicity_high_genus}, the map $f$ is injective. Moreover, if $\phi \in \emcg$ is the lift of the mapping class of a homeomorphism of $\Sigma_g$ that extends to $H_g$ through $\alpha \colon \Sigma_g \rightarrow \partial H_g$, then $ \rho_{p,g}(\phi)\cdot v_0 = v_0$. Denote by $\eemcg\subset \emcg$ the subgroup generated by these elements $\phi$. By definition, we have
\begin{gather*} f ( (\mathcal{A}_{p,g})' ) \subset \big\{ v\in V_p(\Sigma_g) \mbox{ such that } \rho_{p,g}(\phi)\cdot v = v, \mbox{ for all }\phi\in \eemcg \big\}.\end{gather*}

In particular, for any trivalent banded graph $\Gamma$, we have the inclusion $ f ( (\mathcal{A}_{p,g})' ) \subset W_{[0]}(\Gamma)$ where $[0]$ is the class of the coloring sending every edge of $\Gamma$ to~$0$.
\begin{proof}[Proof of Theorem~\ref{main_th} when $p=2r^2$ and $p=18$] When $p=2r^2$, with $r$ an odd prime or $p=18$, then $W_{[0]}$ is one-dimensional, generated by $v_0$. Hence we have the equalities $(\mathcal{A}_{18, g})' = \{ \mathds{1} \}$ and $(\mathcal{A}_{2r^2, 2})' = \{ \mathds{1} \}$. Then the Schur lemma proves our claim.
\end{proof}

\subsection[The case when $p=4r$]{The case when $\boldsymbol{p=4r}$}\label{section4.1}

Assume that $p=4r$ with $r$ an odd prime and write $k:=2r-2$. By Lemma \ref{lemmachoucroute}, a color $i\in I_{p}$ satisfies $\mu_i = 1$ if and only if either $i=0$ or $i=k$. Consider a framed link $L\subset \Sigma_g\times \{\frac{1}{2}\}\subset \Sigma_g\times [0,1]$, thickened inside $\Sigma_g\times \{\frac{1}{2} \}$, and color $L^p $ by $\omega$. In \cite[Section~7]{BHMV2}, it is shown that the operators $\Add_p(L^p(\omega))$ and $\Add_p(L(S_k))$ are equal and only depend on the homology class of $L$ in $\HSg$. Hence we have an injective morphism of algebras
\begin{gather*} i \colon \ \cHSg \hookrightarrow \mathcal{A}_{p,g}.\end{gather*}
The action of a mapping class in homology induces a surjective group morphism $p\colon \emcg \rightarrow \Spg$. By definition of the Witten--Reshetikhin--Turaev representations, we have the following Egorov identity:
\begin{gather}\label{Egorov} \rho_{p,g}(\phi)^{-1}i(w)\rho_{p,g}(\phi) = i ( p(\phi)\cdot w ), \qquad\! \text{for all} \quad\! \phi\in \emcg,\quad\! w\in \cHSg. \!\!\! \end{gather}

Consider a genus $g$ banded trivalent graph $\Gamma$ obtained by connecting $g$ trivial framed knots by a trivalent banded tree. Such a graph is called a \textit{Lollipop graph}, the $g$ edges corresponding to the trivial framed links are called \textit{loop edges} and the edges of the tree are called \textit{trunk edges}. Since $(k,k,k)$ is not $p$-admissible, the space $W_{[0]}(\Gamma)$ is spanned by the vectors $u_{\sigma}$ associated to colorings $\sigma$ such that $\sigma(e)=0$ if $e$ is a trunk edge and $\sigma(e)\in \{0,k \}$ if $e$ is a loop edge. This basis is in natural bijection with the elements of $H_1(H_g, \mathbb{Z}/2\mathbb{Z})$, thus we have a natural isomorphism $\phi\colon \cHHg \cong W_{[0]}(\Gamma)$. Moreover, the isomorphism $\phi$ is equivariant for the actions of $\cHSg$. Precisely, the following diagram commutes
\begin{equation*} \begin{tikzcd}
\mathbb{C}[H_1(\Sigma_g, \mathbb{Z}/2\mathbb{Z})]\times \mathbb{C}[H_1(H_g, \mathbb{Z}/2\mathbb{Z})]
\arrow[r," "] \arrow[d, "i\times \phi", "\cong"']
& \mathbb{C}[H_1(H_g, \mathbb{Z}/2\mathbb{Z})]
\arrow[d, "\phi","\cong"']
\\ i\left(\mathbb{C}[H_1(\Sigma_g, \mathbb{Z}/2\mathbb{Z})]\right)\times W_{[0]}(\Gamma)
\arrow[r, " "]
& W_{[0]}(\Gamma).
\end{tikzcd}\end{equation*}

We denote by $P$ the orthogonal projector of $V_p(\Sigma_g)$ on the subspace of vectors fixed all operators in $i (\cHSg )$. Clearly $P$ belongs to $(\mathcal{A}_{p,g})'$.

Consider a symplectic basis $(x_i,y_i)_{i=1, \ldots, g}$ of $H_1(\Sigma_g, \mathbb{Z}/2\mathbb{Z})$, that is classes of curves, still denoted $x_i$, $y_i$, such that the algebraic intersection of $x_i$ with $y_j$ is equal to the Kronecker symbol $\delta_{i,j}$ and such that the intersections of $x_i$ and $y_i$ with $x_j$ and $y_j$ are null when $i\neq j$. We also suppose that the homeomorphism $\alpha\colon \Sigma_g\cong \partial H_g$, used to define the spaces $V_p(\Sigma_g)$ in Section~\ref{section2.2}, sends the curves $x_i$ to contractible curves in $H_g$. Define
\begin{gather*}\Theta_i := \tfrac{1}{2} ( -1 + x_i + y_i +x_i y_i ) \in \cHSg. \end{gather*}
The $\Theta_i$'s pairwise commute, satisfy $\Theta_i^2=1$ and $i(\Theta_i)\cdot v_0$ is the vector $u_{\sigma}$ associated to the coloring $\sigma$ sending the $i$-th loop edge to $k$ and other edges to $0$. In particular, writing $\mathcal{W}_g:= \mathbb{C}[\Theta_1, \ldots, \Theta_g]\subset \mathbb{C}[H_1(\Sigma_g,\mathbb{Z}/2\mathbb{Z})]$, we have $i(\mathcal{W}_g) \cdot v_0 = W_{[0]}(\Gamma)$.

Denote by $I\subset \cHSg$ the ideal generated by the elements $(x_i-1)$ for $1\leq i \leq g$.
\begin{Lemma}\label{lemma_4r}Consider the action of $\Spg$ on $i( \cHSg )$. Then:
\begin{enumerate}\itemsep=0pt
\item[$1.$] The subspace of vectors fixed by $\Spg$ is $\Span (\mathds{1}, P)$.
\item[$2.$] For every $w\in \mathcal{W}_g$ and $\phi \in \Spg$, we have $ \phi\cdot w - w \in I $.
\end{enumerate}
\end{Lemma}
\begin{proof}The first point follows from the well-known fact that the action of $\Spg$ on $\HSg$ has two orbits: the singleton containing the neutral element and the set containing the other elements. To prove the second point, denote by $X_i$, $Y_i$, $Z_{i,j}$ for $1\leq i,j \leq g$ the classes in $\HSg$ of the Dehn twists of Fig.~\ref{gen_twists} generating $\Spg$. We suppose that $X_i$ and $Y_i$ represent the classes of Dehn twists along curves whose homology classes are $x_i$ and $y_i$ respectively oriented such that $X_i\cdot y_i =x_i y_i$ and $Y_i\cdot x_i=x_iy_i$. First note that the elements $\Theta_i$ are invariant under the action of the $X_j$ and~$Y_j$. We are reduced to show that for any $w\in \mathcal{W}_g$, we have $Z_{i,j}\cdot w - w \in I$. We make the proof for the generator $Z_{1,2}$, the other cases are similar.

\begin{figure}[!h]\centering
\includegraphics[width=8.0cm]{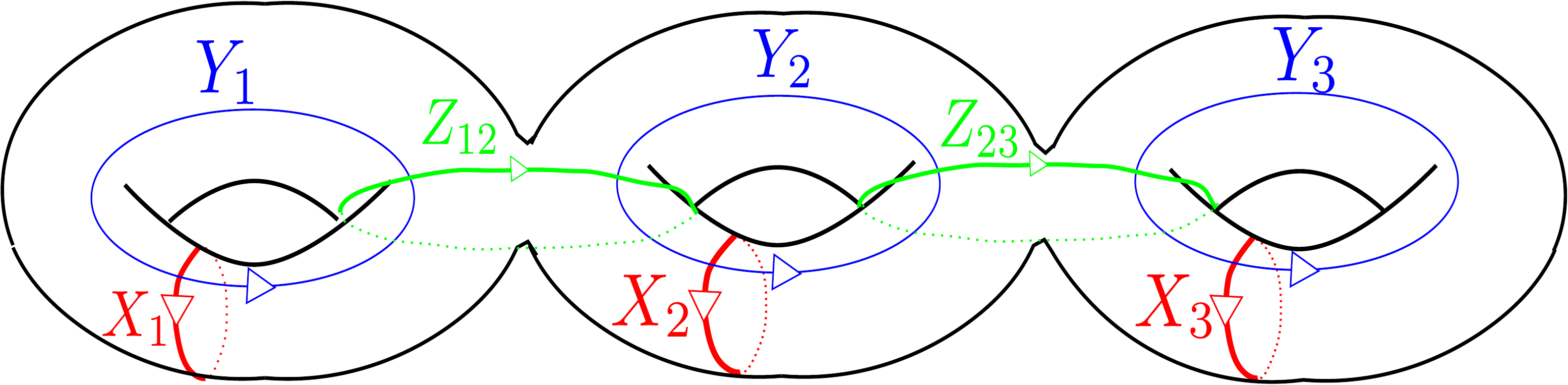}
\caption{The oriented curves defining some Dehn twists whose homology classes generate $\operatorname{Sp}(6, \mathbb{Z}/2\mathbb{Z})$.}\label{gen_twists}
\end{figure}

First note that $Z_{1,2}\cdot \Theta_i = \Theta_i$ when $i\notin\{1,2\}$. Then we compute:
\begin{gather*}
Z_{1,2}\cdot \Theta_1 -\Theta_1 = \tfrac{1}{2}y_1(1+x_1)(x_1x_2-1) \in I, \\
Z_{1,2}\cdot \Theta_2 -\Theta_2 = \tfrac{1}{2}y_2(1+x_2)(x_1x_2-1) \in I.\tag*{\qed}
\end{gather*}\renewcommand{\qed}{}
\end{proof}

The case $p=4r$ of Theorem~\ref{main_th} follows from the
\begin{Proposition}\label{prop_case_4r}If $p=4r$, with $r$ an odd prime, is such that $v_0 \in V_p(\Sigma_g)$ is cyclic, then $V_p(\Sigma_g)$ is the direct sum of two simple submodules.
\end{Proposition}

\begin{proof}Consider the linear map $h\colon \mathcal{W}_g \rightarrow W_{[0]}(\Gamma)$ defined by $h(w)=i(w)\cdot v_0$. The map $h$ is surjective and its kernel is the ideal $I$. Lemma~\ref{lemma_4r}(2) and equation~\eqref{Egorov} imply that
\begin{gather*} \rho_{p,g}(\phi)\circ i(w)\cdot v_0 = i(w)\circ \rho_{p,g}(\phi)\cdot v_0,\qquad \text{for all} \quad \phi \in \emcg, \quad w\in \mathcal{W}_g.\end{gather*}
Let $\theta\in (\mathcal{A}_{p,g})'$. Since $\theta\cdot v_0$ lies in $W_{[0]}(\Gamma)$ and $h$ is surjective, there exists an element $w\in \mathcal{W}_g$ such that $i(w)\cdot v_0 = \theta \cdot v_0$. Moreover, if $\phi \in \emcg$, then
\begin{gather*} \theta \circ \rho_{p,g}(\phi) \cdot v_0 = \rho_{p,g}(\phi)\circ \theta \cdot v_0 = \rho_{p,g}(\phi)\circ i(w) \cdot v_0 = i(w)\circ \rho_{p,g}(\phi)\cdot v_0.\end{gather*}
The cyclicity of $v_0$ implies that $\theta = i(w) \in i(\mathcal{W}_g)$, hence we have the inclusion $\left(\mathcal{A}_{p,g}\right)'\subset i(\mathcal{W}_g)\subset i\left( \mathbb{C}[H_1(\Sigma_g, \mathbb{Z}/2\mathbb{Z})]\right)$.

Moreover, Lemma \ref{lemma_4r}(1) implies that $(\mathcal{A}_{p,g})' \cap i( \mathbb{C}[H_1(\Sigma_g, \mathbb{Z}/2\mathbb{Z})]) = \Span(\mathds{1}, P)$. Hence we have the equality $(\mathcal{A}_{p,g})'=\Span(\mathds{1},P)$ which proves the claim.
\end{proof}

\subsection[The case when $p=2r_1r_2$]{The case when $\boldsymbol{p=2r_1r_2}$}

Assume that $p=2r_1r_2$ with $r_1$, $r_2$ distinct odd primes. In this case, there exists an unique integer $x\in \{ 1, \ldots, r_1r_2-2\}$ such that $\muu{x}=1$. Then $x$ is even and verifies either
\begin{gather*}\begin{cases}
x\equiv -2 &\pmod{r_1}, \\
x\equiv 0 & \pmod{r_2},
\end{cases} \qquad \text{or}\qquad \begin{cases}
x\equiv 0 &\pmod{r_1}, \\
x\equiv -2 & \pmod{r_2}.
\end{cases}\end{gather*}
We begin by stating a technical lemma whose proof will be the subject of the next subsection:
\begin{Lemma}\label{GA}
 If $(x,x,x)$ is $p$-admissible and $r_1,r_2>37$, then we have the following inequality
\begin{gather*} \sixj{x}{x}{2}{x}{x}{0} \sixj{x}{x}{4}{x}{x}{x} \neq \sixj{x}{x}{4}{x}{x}{0}\sixj{x}{x}{2}{x}{x}{x}.\end{gather*}
\end{Lemma}

 Consider the two banded graphs $\Gamma_1= \Thetagraph$ and $\Gamma_2=\Eyesgraph$. Denote by $b_1, b_2 \in V_p(\Sigma_2)$ the two vectors $b_1:= \eyesgraph{$x$}{0}{0}$ and $b_2 := \eyesgraph{0}{0}{$x$}$.

\begin{Lemma}\label{lemma_2r1r2}
Suppose $p=2r_1r_2$, with $r_1$, $r_2$ distinct odd primes, such that either $(x,x,x)$ is not $p$-admissible, or $r_1,r_2>37$, then
\begin{gather*} W_{[0]}(\Gamma_1) \bigcap W_{[0]}(\Gamma_2) \subset \Span ( v_0, b_1, b_2 ).\end{gather*}
\end{Lemma}

\begin{proof}The subspaces $W_{[0]}(\Gamma)$ are spanned by the vectors associated to colorings of the edges of~$\Gamma$ by the elements $0$ and $x$. If $(x,x,x)$ is not $p$-admissible, then these colorings of $\Gamma_2$ correspond to the elements $v_0$, $b_1$, $b_2$ and the proof is immediate. Suppose $(x,x,x)$ is $p$-admissible. We need to show that if $v=\alpha \thetagraph{$x$}{0}{$x$} + \beta \thetagraph{$x$}{$x$}{$x$} \in W_{[0]}(\Gamma_2)$, then $\alpha=\beta = 0$.

Lemma \ref{fusion_rule} implies that
\begin{gather*}
v= \left( \alpha \sixj{x}{x}{2}{x}{x}{0} +\beta \sixj{x}{x}{2}{x}{x}{x} \right) \eyesgraph{x}{2}{x}\\
\hphantom{v=}{}+ \left( \alpha \sixj{x}{x}{4}{x}{x}{0} +\beta \sixj{x}{x}{4}{x}{x}{x} \right) \eyesgraph{x}{4}{x} + v',
\end{gather*}
where $v'$ is a linear combination of vectors of the form $\eyesgraph{$x$}{$i$}{$x$}$ for $i\neq 2,4$. Since $v\in W_{[0]}(\Gamma_2)$, the two coefficients in front of $\eyesgraph{$x$}{2}{$x$}$ and $\eyesgraph{$x$}{4}{$x$}$ vanish. Hence we get
\begin{gather*}\begin{pmatrix} \sixj{x}{x}{2}{x}{x}{0} & \sixj{x}{x}{2}{x}{x}{x} \\ \sixj{x}{x}{4}{x}{x}{0} & \sixj{x}{x}{4}{x}{x}{x} \end{pmatrix} \cdot \begin{pmatrix} \alpha \\ \beta \end{pmatrix} = \begin{pmatrix} 0\\0 \end{pmatrix}.\end{gather*}

Lemma \ref{GA} concludes the proof.
\end{proof}

\begin{Lemma}\label{lemmalacon}
There exists an element $a \in \mathcal{A}_{2r_1r_2,1}$ such that
\begin{gather*}
a\cdot u_0 = u_x, \qquad a \cdot u_x = u_0.
\end{gather*}
\end{Lemma}

\begin{proof} It is enough to show that there exists an operator $\psi \in (\mathcal{A}_{2r_1r_2,1})'$ such that
\begin{gather*} \psi \cdot u_0 = u_x \qquad \text{and} \qquad \psi \cdot u_x =u_0.\end{gather*}

The cyclicity of $u_0$, provided by Lemma \ref{lemma_genusone1}, implies the existence of $a\in\mathcal{A}_{2r_1r_2,1}$ such that
$a\cdot u_0 = u_x$. If such a $\psi$ exists, then
\begin{gather*} a\cdot u_x = a \circ \psi \cdot u_0 = \psi \circ a \cdot u_0 = u_0.\end{gather*}

The operator $\psi$ is defined as follows. Given $i\in \{0, \ldots, r_1r_2-2\}$, only one of the following two cases occurs:
\begin{enumerate}\itemsep=0pt
\item Either there exists $j\in\{0, \ldots, r_1r_2-2 \}$ so that
\begin{gather*} \begin{cases}
j\equiv i & \pmod{2r_1}, \\
j\equiv -i-2 & \pmod{r_2},
\end{cases}\end{gather*}
 and we set $\psi(u_i):=u_j$.
\item Or there exists $j\in\{0, \ldots, r_1r_2-2 \}$ so that
\begin{gather*} \begin{cases}
j\equiv i & \pmod{2r_2}, \\
j\equiv -i-2& \pmod{r_1},
\end{cases}\end{gather*}
 and we set $\psi(u_i):=-u_j$.
\end{enumerate}
A straightforward computation shows that $\psi$ commutes with $\rho_{p,1}(T)$ and $\rho_{p,1}(S)$ and either $\psi$ or $-\psi$ sends $u_0$ to $u_x$.
\end{proof}

\begin{figure}[!h]\centering
\includegraphics[width=12cm]{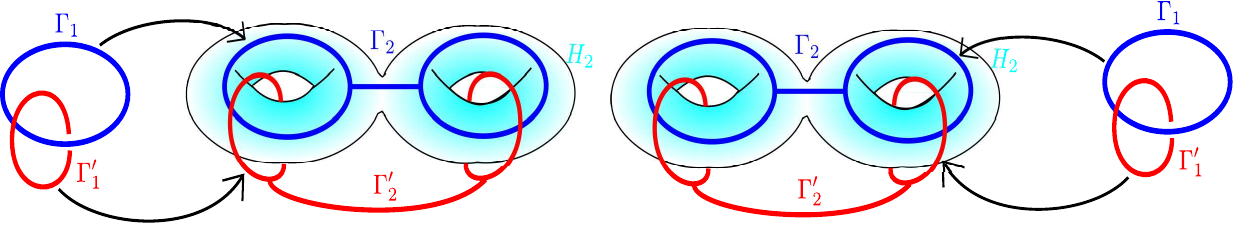}
\caption{The genus $2$ graphs $\Gamma_2$ and $\Gamma'_2$ and the genus one subgraphs $\Gamma_1$ and $\Gamma'_1$ used to define $a\otimes \mathds{1}$ on the left and $\mathds{1}\otimes a$ on the right.}\label{graph_embedding}
\end{figure}

We define two operators $a\otimes \mathds{1}, \mathds{1}\otimes a \in \mathcal{A}_{p, 2}$ as follows. Consider $\Gamma_2, \Gamma'_2 \subset S^3$ the two entangled genus $2$ banded graphs of Fig.~\ref{graph_embedding}. Denote by $\Gamma_1\subset \Gamma_2$ and $\Gamma'_1\subset \Gamma'_2$ the genus one subgraphs of the left picture of Fig.~\ref{graph_embedding}. Fix $H_2, H'_2$ some tubular neighborhood of $\Gamma_2$ and $\Gamma'_2$ respectively which do not intersect and $H_1\subset H_2$, $H'_1\subset H'_2$ some tubular neighborhoods of $\Gamma_1$ and $\Gamma'_1$ respectively. Identify the closure of $S^3\setminus (H_i\cup H'_i)$ with $\Sigma_i\times [0,1]$ and $V_p(\Sigma_i)$ with the space of linear combinations of framed links in $H_i$ quotiented by the kernel of the Hopf pairing induced by the Heegaard splitting $S^3= H_i \cup \overline{(S^3\setminus H_i)}$. From Lemma~\ref{lemmalacon}, there exists some $a\in \mathcal{A}_{p,1}$ such that $a\cdot u_x = u_0$ and $a\cdot u_0=u_x$. Let $\mathcal{L}$ be a linear combination of aligned links in $\overline{S^3\setminus (H_1\cup H'_1)}\cong \Sigma_1\times [0,1]$ such that $\Add_{p,1}\left(\mathcal{L}(\omega)\right)=a$. Composing eventually by an isotopy, we can suppose that the framed parallel links of $\mathcal{L}$ are in $\overline{S^3\setminus (H_2\cup H'_2)}\cong \Sigma_2\times [0,1]$. We denote by $a\otimes \mathds{1}:=\Add_{p,2}\left(\mathcal{L}(\omega)\right) \in \mathcal{A}_{p,2}$ the resulting operator. Note that this definition is not canonical since it depends on the choice of $a$ and on how the middle edges of~$\Gamma_2$ and~$\Gamma'_2$ are entangled with the framed links of $\mathcal{L}$. Nevertheless, by definition this operator satisfies $(a\otimes \mathds{1})\cdot v_0 = b_1$ and $(a\otimes \mathds{1})^2 \cdot v_0 = v_0$. Similarly, by considering the genus one subgraphs of $\Gamma_2$ and $\Gamma'_2$ of the right picture of Fig.~\ref{graph_embedding}, we define an operator $\mathds{1}\otimes a \in \mathcal{A}_{p,2}$ such that $(\mathds{1}\otimes a)\cdot v_0 = b_2$ and $(\mathds{1}\otimes a)^2 \cdot v_0 = v_0$. We further suppose that the aligned framed links defining these two operators are not entangled so the two operators $a\otimes \mathds{1}$ and $\mathds{1}\otimes a$ commute. We eventually define $a\otimes a := (a\otimes \mathds{1})\circ (\mathds{1}\otimes a)$ which satisfies $(a\otimes a)\cdot v_0 = \eyesgraph{$x$}{0}{$x$}$ and $(a \otimes a)^2\cdot v_0 = v_0$.

\begin{proof}[Proof of Theorem \ref{main_th} when $p=2r_1r_2$]
Recall we defined a linear map $f\colon (\mathcal{A}_{p,2})' \hookrightarrow V_p(\Sigma_2)$ by $f(\theta)=\theta\cdot v_0$ and that the image $f\left( (\mathcal{A}_{p,2})' \right) \subset W_{[0]}(\Gamma_2)$ is invariant under the orientation-preserving homeomorphisms of the handlebody $H_2$. Using Lemma~\ref{lemma_2r1r2} and the fact that the vectors in the image of $f$ must be invariant under a homeomorphism permuting the two handles of $H_2$, we have the inclusion
\begin{gather*} f ( (\mathcal{A}_{p,2})' ) \subset \Span ( v_0, b_1+b_2 ).\end{gather*}

We need to show that $f( (\mathcal{A}_{p,2})' )= \mathbb{C}\cdot v_0$. By contradiction, suppose that there exists $\theta \in (\mathcal{A}_{p,2})'$ such that $\theta \cdot v_0 = b_1+b_2= (a\otimes \mathds{1} + \mathds{1}\otimes a) \cdot v_0$. Then we have
\begin{gather*}
\theta ^2 \cdot v_0 =( a \otimes \mathds{1} + \mathds{1}\otimes a )^2\cdot v_0= 2 v_0 + 2\eyesgraph{$x$}{0}{$x$}.
\end{gather*}
Thus $\theta^2\cdot v_0$ does not belong to $ \Span ( v_0, b_1+b_2 ) $. This contradicts the fact that $\theta^2 \in (\mathcal{A}_{p,2})'$.
\end{proof}

\subsection{Proof of Lemma \ref{GA}}

In this subsection we consider $p=2r_1r_2$, with $r_1, r_2>37$ two distinct odd primes. We suppose that there exists $x\in\{1,\ldots, r_1r_2-2\}$ so that $(x,x,x)$ is $p$-admissible and
\begin{gather*}\begin{cases}
x\equiv 0 & \pmod{2r_1}, \\
x\equiv -2 & \pmod{r_2}.
\end{cases}
\end{gather*}

We also choose some primitive $r_1$-{th} and $r_2$-{th} roots of unity $A_1$ and $A_2$, such that $A^2=A_1A_2$. In particular, we have $A^{2x}=A_2^{-2}$.

The goal of this subsection is to show that
\begin{gather*} D:= \sixj{x}{x}{2}{x}{x}{0}\sixj{x}{x}{4}{x}{x}{x} - \sixj{x}{x}{4}{x}{x}{0}\sixj{x}{x}{2}{x}{x}{x} \neq 0.\end{gather*}

A straightforward computation, using the formula of the recoupling coefficients~\cite{MV}, gives
\begin{gather*}
D = (-1)^{\frac{x}{2}+1}\frac{[3][5]![x] \left[\frac{3}{2}x+1 \right]!(\left[\frac{x}{2}\right]!)^3}{[2][x+3]!([x+2]!)^2[x+3] \left[\frac{x}{2}+1\right]}
\left( \left[\frac{x}{2}-1\right]^2 \left[\frac{x}{2}\right]^2\left[\frac{x}{2}+1\right] \right.\\
\hphantom{D =}{} - [2]^2\left[\frac{x}{2}\right]^3\left[\frac{3}{2}+2\right] \left[\frac{x}{2}+1\right] + \left[\frac{x}{2}-1\right]\left[\frac{x}{2}\right] \left[\frac{x}{2}+1\right] \left[\frac{x}{2}+2\right] \left[\frac{x}{2}+3\right]\\
\left.\hphantom{D =}{} + [x-1][x+3] \left[\frac{x}{2}\right]^2\left[\frac{x}{2}+1\right] - [x-1]\left[\frac{3}{2}x+2\right][x+3] \right) \\
\hphantom{D}{} = (-1)^{\frac{x}{2}+1}\frac{[3][5]![x]\left[\frac{3}{2}x+1\right]!(\left[\frac{x}{2}\right]!)^3}{[2][x+3]!([x+2]!)^2[x+3]\left[\frac{x}{2}+1\right](A^2-A^{-2})^7A_1^{10}A_2^{9}}\cdot P(A_1,A_2),
\end{gather*}
where $P$ is defined by
\begin{gather*}
P(x,y):= x^{20}y^{16}-x^{17}y^{17}-x^{16}y^{18}+x^{19}y^{13}-4x^{18}y^{14}+3x^{17}y^{15}+2x^{15}y^{17}-x^{19}y^{11}\\
\hphantom{P(x,y):=}{} -5x^{17}y^{13}-4x^{15}y^{15}+4x^{14}y^{16}-2x^{13}y^{17}-x^{18}y^{10} +2x^{17}y^{11}+6x^{16}y^{12} \\
\hphantom{P(x,y):=}{}+2x^{15}y^{13}-x^{14}y^{14}+2x^{13}y^{15}+x^{11}y^{17}+2x^{15}y^{11} +x^{14}y^{12}+x^{13}y^{13}-6x^{12}y^{14} \\
\hphantom{P(x,y):=}{}+x^{10}y^{16}+x^{17}y^7+4x^{16}y^8-x^{15}y^9-4x^{14}y^{10}+4x^{12}y^{12}+x^{15}y^7+x^{13}y^9 \\
\hphantom{P(x,y):=}{}-4x^{12}y^{10}+x^{11}y^{11}+4x^{10}y^{12}-4x^8y^{14}-x^7y^{15}-2x^{15}y^5-6x^{14}y^6-4x^{13}y^7 \\
\hphantom{P(x,y):=}{}+2x^{12}y^8 -8x^{11}y^9-6x^{10}y^{10}-6x^9y^{11}-x^7y^{13} +x^{13}y^5+6x^{11}y^7+6x^{10}y^8 \\
\hphantom{P(x,y):=}{}+8x^9y^9-2x^8y^{10}+4x^7y^{11}+6x^6y^{12}+2x^5y^{13}+x^{13}y^3+4x^{12}y^4 -4x^{10}y^6 \\
\hphantom{P(x,y):=}{}-x^9y^7+4x^8y^8 -x^7y^9-x^5y^{11}-4x^8y^6+4x^6y^8+x^5y^9-4x^4y^{10}-x^3y^{11} \\
\hphantom{P(x,y):=}{}-x^{10}y^2+6x^8y^4-x^7y^5-x^6y^6-2x^5y^7-x^9y-2x^7y^3 +x^6y^4-2x^5y^5\\
\hphantom{P(x,y):=}{}-6x^4y^6-2x^3y^7+x^2y^8+2x^7y-4x^6y^2+4x^5y^3+5x^3y^5+xy^7-2x^5y\\
\hphantom{P(x,y):=}{}-3x^3y^3+4x^2y^4-xy^5+x^4+x^3y-y^2.
\end{gather*}

Note that $P$ does not depend neither on $r_1$, $r_2$ nor on $x$. We have to show that $P(A_1,A_2)\neq 0$. By contradiction, suppose that $P(A_1, A_2)=0$ and define $H(X):=P(X^{r_2}, X^{r_1})\in \mathbb{Z}[X]$. Since~$P$ has integral coefficients the formula $P(A_1,A_2)=0$ holds for any roots of unity $A_1$ and $A_2$ of order $r_1$ and $r_2$ respectively. In particular, we have $H(A^2)=0$ and $H(X)$ is divisible by the product $\phi_{r_1}(X)\phi_{r_2}(X)$ of the two cyclotomic polynomials of order $r_1$ and $r_2$. Thus its degree satisfies $\deg (H) \geq (r_1-1)(r_2-1)$. Using the above formula for $P$ we deduce that
\begin{gather*}
 (r_1-1)(r_2-1) \leq \deg (H) = \begin{cases}
20r_2+16r_1, & \text{if} \ 2r_2\geq r_1, \\
16r_2 +18r_1, &\text{if} \ 2r_2\leq r_1,
\end{cases}\\
\qquad {} \Leftrightarrow
 \begin{cases}
(r_1-21)(r_2-17)\leq 356, & \text{if} \ 2r_2\geq r_1, \\
(r_1-17)(r_2-19)\leq 352, &\text{if} \ 2r_2\leq r_1.
\end{cases}
\end{gather*}

Since $r_1, r_2> 37$ are prime, then $r_1, r_2\geq 41$. We conclude by noticing that none of the two equations of the above system does have solutions satisfying $r_1, r_2\geq 41$.

\section{Partial generalizations in higher genus}

The main obstruction to extend the above techniques in genus $g\geq 3$ lies in the proof of Proposition \ref{prop_cyclicity_high_genus}. To prove the cyclicity of the vacuum vector, we found linear independent equations satisfied by any vector of $W_{[\sigma]}\bigcap \left( A_{p,g}\cdot v_0\right)^{\bot} $ and deduced that this space must be trivial. The non-triviality of these equations followed from the fact that the Hopf pairing $( u_{\sigma}, v_0)^H_{p,2}$ is a $3j$-symbol, thus does not vanish. In higher genus, for a suitable choice of trivalent graph, this Hopf pairing is a product of $6j$-symbols, which might a priori vanish. Indeed, suppose $p=2(k+2)$ with $k\geq 2$ an even integer, and consider the two following families of $6j$-symbols:
$$
\begin{tabular}{|c|c|}
\hline
\text{Type I} & \text{Type II}
\\ \hline
 $\left< \Tetahedrongraph{$\frac{k}{2}$}{$\frac{k}{2}$}{$\frac{k}{2}$}{$a$}{$b$}{$c$} \right> $
 &
$ \left< \Tetahedrongraph{$\frac{k}{2}$}{$a$}{$k-a$}{$b$}{$c$}{$b$} \right>$
 \\
 \text{ where }$a+b+c\equiv 2 \pmod{4}$ & \text{ where }$a+\frac{c+k}{2} \equiv 1 \pmod{2}$
\\ \hline
\end{tabular}
$$
\begin{Proposition}\label{prop_null_6j}
The $6j$-symbols of type I and type II, defined above, vanish.
\end{Proposition}

We postpone the technical proof of Proposition \ref{prop_null_6j} to Appendix~\ref{appendixA}. We conjecture that the $6j$-symbols of type I and II are the only vanishing $6j$-symbols:
\begin{Conjecture}\label{Conjecture}If $4$ divides $p$, then the only vanishing $6j$-symbols at level $p$ are the two families type~I and type~II from above. If~$4$ does not divide~$p$, there is no vanishing $6j$-symbol at level~$p$.
\end{Conjecture}

Numerical computations permitted us to verify the validity of the conjecture up to level $66$. Motivated by this conjecture, we will show the following:

\begin{Proposition}\label{prop_high_genus} \quad
\begin{enumerate}\itemsep=0pt
\item[$1.$] Let $p=2r_1r_2$ with $r_1$, $r_2$ distinct odd primes such that either $r_1,r_2>37$ or the element $x$ defined in the introduction satisfies $3x>2r_1r_2-2$. Suppose that $p$ verifies Conjecture~{\rm \ref{Conjecture}}. Then $V_p(\Sigma_g)$ is simple for arbitrary $g\geq 3$.
\item[$2.$] Let $p=4r$ with $r$ an odd prime and suppose that $p$ verifies Conjecture~{\rm \ref{Conjecture}}. Then $V_p(\Sigma_3)$ is the direct sum of two simple modules.
\end{enumerate}
\end{Proposition}

Theorem \ref{theorem_high_genus} is a consequence of this proposition and the computer assisted verification of Conjecture~\ref{Conjecture} up to level~$66$. To prove Proposition~\ref{prop_high_genus}, we first show that the arguments of section $4$ generalize:
\begin{Lemma}\label{lemma_highgen1}
Suppose that $p, g \geq 3$ are such that $v_0\in V_p(\Sigma_g)$ is cyclic, then:
\begin{enumerate}\itemsep=0pt
\item[$1.$] If $p=2r_1r_2$ with $r_1$, $r_2$ distinct odd primes such that either $r_1,r_2>37$ or $3x>2r_1r_2-4$, then $V_p(\Sigma_g)$ is simple.
\item[$2.$] If $p=4r$ with $r$ an odd prime, then $V_{4r}(\Sigma_g)$ is the sum of two simple modules.
\end{enumerate}
\end{Lemma}

\begin{proof}The case $p=4r$ follows from Proposition~\ref{prop_case_4r}. Suppose that $p=2r_1r_2$ is such that $v_0\in V_p(\Sigma_g)$ is cyclic and recall that there exists an unique color $x\neq 0$ such that $\mu_x=1$.

\textbf{Step 1:} Suppose that $(x,x,x)$ is $p$-admissible and $r_1,r_2>37$. Let $\Gamma_1$, $\Gamma_2$ be two trivalent graphs which only differ by a single Whitehead move inside a ball~$B^3$, as drawn in Fig.~\ref{whitehead_move}. Then
\begin{gather*} W_{[0]}(\Gamma_1)\cap W_{[0]}(\Gamma_2) \subset \Span \big( u_{\sigma}^{\Gamma_1}, \ \text{such that} \ \sigma(a)\sigma(b)\sigma(c)\sigma(d)=0 \big).\end{gather*}
Indeed, let $\sigma_1$, $\sigma_2$ be two $p$-admissible colorings of $\Gamma_1$, with colors~$0$ or~$x$, such that
$ \sigma_1(e) = \sigma_2(e)$, for all $e \in E(\Gamma_1)\setminus \{i \} $ and with $\sigma_i(a)=\sigma_i(b)=\sigma_i(c)=\sigma_i(d)=x$ and $\sigma_1(i)=0$, $\sigma_2(i)=x$.

Suppose that there exists $(\alpha, \beta) \in \mathbb{C}^2$ such that the vector $ v:= \alpha u_{\sigma_1} + \beta u_{\sigma_2} \in W_{[0]}(\Gamma_2)$. Then, as in the proof of Lemma~\ref{lemma_2r1r2}, using fusion-rules, we show that $(\alpha, \beta)$ satisfies a system of two linear equations which are linearly independent by Lemma~\ref{GA}. Hence we have $v=0$.

\textbf{Step 2:} If $i\in \{1,\ldots, g\}$, we denote by $b_i\in V_p(\Sigma_g)$ the vector representing a single framed knot colored by $x$ around the $i$-{th} hole. We have:
\begin{gather*} \bigcap_{\Gamma } W_{[0]}(\Gamma) = \Span (v_0, b_i, 1\leq i \leq g ),\end{gather*}
where the intersection is over the set of embedded banded trivalent graphs of genus~$g$.

Indeed, let $\Gamma$ represents the graph of Fig.~\ref{lunettesgraph} and $\sigma$ a $p$-admissible coloring such that:
\begin{enumerate}\itemsep=0pt
\item[1)] $\sigma(e) \in \{ 0, x \}$, for all $e\in E(\Gamma)$,
\item[2)] 
there exists $i<j$ with $\sigma(a_i)=\sigma(b_i)=\sigma(a_j)=\sigma(b_j)=x$.
\end{enumerate}

We can suppose that for every $i<k<j$, we have $\sigma(a_k)\sigma(b_k)=0$.

Step $1$ with $a=a_i$, $b=a_j$, $c=b_i$ and $d=b_j$ implies that the projection of $u_{\sigma}$ on $\bigcap_{\Gamma} W_{[0]}(\Gamma)$ vanishes. We conclude by noticing that if $\sigma$ is a $p$-admissible coloring of $\Gamma$, with colors in $\{0, x\}$, that does not satisfy the condition~2, then either $u_{\sigma}=b_i$ for some $i\in \{1, \ldots, g\}$ or else $u_{\sigma}=v_0$.

 \begin{figure}[!h]\centering
\includegraphics[width=7.5cm]{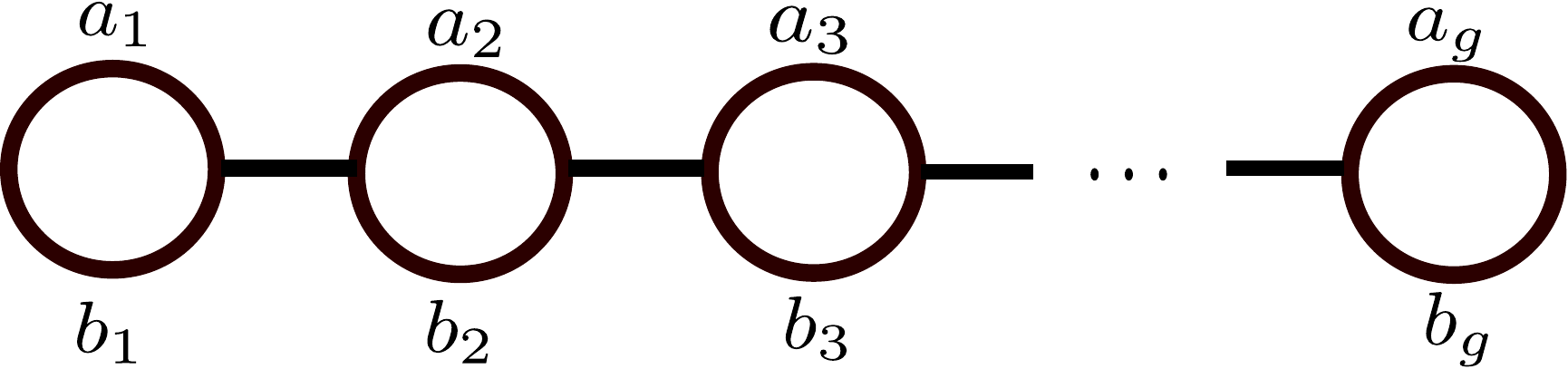}
\caption{A trivalent banded graph of genus $g$. }\label{lunettesgraph}
\end{figure}

\textbf{Step 3:} We conclude the proof as in the case $g=2$. Using Step $2$ and the fact that the vectors of the image of $f$ must be invariant under permutation of the handles, we have the inclusion:
$ f( (\mathcal{A}_{p,g})') \subset \Span ( v_0, b_1+\dots +b_g )$. We must show that $f\left( (\mathcal{A}_{p,g})'\right)=\mathbb{C}\cdot v_0$. By contradiction, suppose that there exists $\theta \in (\mathcal{A}_{p,g})'$ such that
\begin{gather*}
 \theta \cdot v_0 = b_1+\cdots + b_g = ( a \otimes \mathds{1} \otimes \!\cdots\! \otimes \mathds{1}+ \mathds{1}\otimes a \otimes \mathds{1}\otimes \!\cdots\! \otimes \mathds{1} + \dots + \mathds{1}\otimes \!\cdots\! \otimes \mathds{1} \otimes a )\cdot v_0,
\end{gather*}
where $a\otimes \mathds{1}\otimes \ldots \otimes \mathds{1}$ is defined using Lemma~\ref{lemmalacon} as in the $g=2$ case. We have
\begin{gather*}
\theta ^2 \cdot v_0 = \big( a \otimes \mathds{1} \otimes \dots \otimes \mathds{1} + \dots + \mathds{1}\otimes \dots \otimes \mathds{1} \otimes a \big)^2\cdot v_0 \\
\hphantom{\theta ^2 \cdot v_0}{} = g v_0 + 2 ( a\otimes a \otimes \mathds{1}\otimes \dots \otimes \mathds{1} )\cdot v_0 + \dots + 2 ( \mathds{1}\otimes \dots \otimes \mathds{1} \otimes a\otimes a)\cdot v_0.
\end{gather*}
Hence $\theta^2\cdot v_0$ does not belong to $ \bigcap_{\Gamma} W_{[0]}(\Gamma) $. This contradicts the fact that $\theta^2 \in (\mathcal{A}_{p,g})'$.
\end{proof}

\begin{figure}[!h]\centering
\includegraphics[width=7cm]{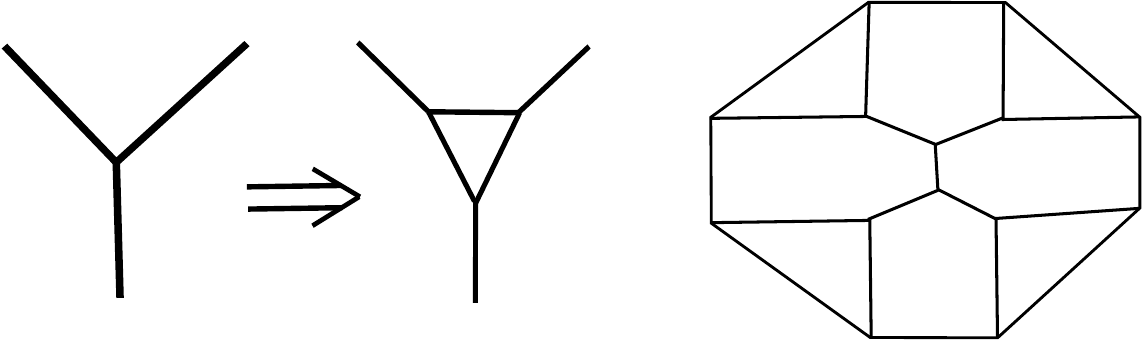}
\caption{On the left, the operation transforming a fly-eyes-graph of genus $g$ into a one of genus $g+1$. On the right, an example of a genus $8$ fly-eyes-graph.}\label{trucgenresup}
\end{figure}

We now extend Proposition \ref{prop_cyclicity_high_genus} to higher genus. A genus $g\geq 2$ trivalent banded graph $\Gamma_g$ is called a \textit{fly-eyes-graph} if it is obtained from $\Gamma_2= \Thetagraph$ by a sequence of $g-2$ operations which consist in choosing an arbitrary vertex and inserting a triangle as drawn in the left part of Fig.~\ref{trucgenresup}. The right part of Fig.~\ref{trucgenresup} provides an example of a genus $8$ fly-eyes-graph. For instance, the Lollipop graphs defined in Section~\ref{section4.1} are not fly-eyes-graphs. There is only one genus $3$ fly-eyes-graph which is the tetrahedron graph. It follows from standard fusion-rules that if $\sigma$ is a $p$-admissible coloring of a fly-eyes-graph, then the Hopf pairing $ ( u_{\sigma}, v_0 )_{p,g}^H$ is a product of $6j$-symbols multiplied by some non-vanishing constant.

Fix $g\geq 3$ and embed a fly-eyes-graph $\Gamma_g$ in $S^3$. Denote by $H_g$ the embedded handlebody $H_g:= S^3 \backslash \mathring{N}(\Gamma_g)$ where $N(\Gamma_g)$ denotes a tubular neighborhood of $\Gamma_g$. For each edge $e\in E(\Gamma_g)$, fix a curve $\gamma_e\subset H_g$ which bounds a disc intersecting $\Gamma_g$ only once along $e$. We construct a map
\begin{gather*} w \colon \ {\mathbb{N}}^{E(\Gamma_g)} \rightarrow V_p(\Sigma_g) \end{gather*}
as follows. To $f\colon E(\Gamma_g)\rightarrow \mathbb{N}$ we associate the class in $V_p(\Sigma_g)$ of the framed link consisting of~$f(e)$ parallel copies of~$\gamma_e$, thickened along $\Sigma_g$, for each edge $e\in E(\Gamma_g)$.

\begin{Lemma}\label{lemma_meuh} If $p=4r$, with $r$ an odd prime, or $p=2r_1r_2$, with $r_1$, $r_2$ distinct odd primes, then $w_f\in \mathcal{A}_{p,g}\cdot v_0$ for all $f\in \{0, 1\}^{E(\Gamma_g)}$.
\end{Lemma}

\begin{proof} We will show the stronger result that if $\sigma$ is a $p$-admissible coloring such that there exists an edge $e$ satisfying $\sigma(e)=1$, then $u_{\sigma}\in \mathcal{A}_{p,g}\cdot v_0$. The conclusion of the lemma will follow. Fix an edge $e\in E(\Gamma_g)$ and denote by $_i\Sigma(\gamma_e)_i$ the surface obtained by cutting~$\Sigma_g$ along~$\gamma_e$ and gluing back two discs with one framed puncture colored by $i$. By the colored splitting~\cite[Theorem~1.4]{BHMV2}, we have a canonical isomorphism
\begin{gather*} V_p(\Sigma_g) \cong \oplus_i V_p(_i\Sigma(\gamma_e)_i).\end{gather*}
Denote by $\Mod(\Sigma_{g-1,2})$ the mapping class group of a genus $g-1$ surface with two boundary components and by $\rho_{p,g}^i \colon \widetilde{\Mod}(\Sigma_{g-1, 2})\rightarrow \GL(V_p(_i\Sigma(\gamma_e)_i))$ the Witten--Reshetikhin--Turaev representation associated to $_i\Sigma(\gamma_e)_i$ (see~\cite{BHMV2} for definitions). There is a natural embedding $\widetilde{\Mod}(\Sigma_{g-1, 2}) \hookrightarrow \emcg$ such that the following diagram commutes for every color~$i$:
\begin{equation*}\begin{tikzcd}
 \widetilde{\Mod}(\Sigma_{g-1, 2}) \arrow[d, hook, ""] \arrow[r, "\rho_{p,g}^i"] &
 \GL(V_p(_i\Sigma(\gamma_e)_i)) \arrow[d, hook, ""] \\
 \emcg \arrow[r, "\rho_{p,g}"]
& \GL(V_p(\Sigma_g)).
\end{tikzcd}\end{equation*}

It follows from the main theorem of \cite{KoberdaSantharoubane17} that for $i=1$, the representation $\rho_{p,g}^1$ is irreducible. Lemmas \ref{lemma_genusone1} and \ref{lemma_reduction} imply that there exists a $p$-admissible coloring $\sigma$ with $\sigma(e)=1$ and \mbox{$u_{\sigma} \in \mathcal{A}_{p,g}\cdot v_0$}. The irreducibility of $\rho_{p,g}^1$ implies that the cyclic space of $u_{\sigma}$ contains the whole subspace $V_p(_1\Sigma(\gamma_e)_1)$. This subspace is spanned by the vectors $u_{\sigma}$ such that $\sigma(e)=1$. This concludes the proof.
\end{proof}

Suppose that $p=4r$ and consider the genus $3$ tetrahedron graph~$\Gamma$. A $p$-admissible coloring~$\sigma$ of~$\Gamma$ is said of type I or II if the corresponding $6j$-symbol $\left< u_{\sigma}\right>$ is of type I or II respectively. Denote by $Z_{p,3}\subset V_p(\Sigma_3)$ the subspace generated by the vectors $u_{\sigma}$ associated to $p$-admissible colorings of the tetrahedron graph which are neither of type I nor of type II.
\begin{Lemma}\label{lemma_ultimate}\quad
\begin{enumerate}\itemsep=0pt
\item[$1.$] If $g\geq 3$, $p=2r_1r_2$ with $r_1$, $r_2$ distinct odd primes, and if $p$ verifies Conjecture~{\rm \ref{Conjecture}}, then $v_0\in V_p(\Sigma_g)$ is cyclic.
\item[$2.$] If $p=4r$, with $r$ an odd prime, verifies Conjecture~{\rm \ref{Conjecture}}, then $Z_{p,3}\subset \mathcal{A}_{p,3}\cdot v_0$.
\end{enumerate}
\end{Lemma}

\begin{proof} The proof is similar to the genus $2$ case. Choose a vector
\begin{gather*}
v \in W_{[\sigma]}\bigcap ( A_{p,g}\cdot v_0)^{\bot} v = \sum_{\tau\in [\sigma]}{\alpha_{\tau}u_{\tau}}.
\end{gather*}
If $p=4r$, we further suppose that $v\in Z_{p,3}$. We must find $\dim \left(W_{[\sigma]}\right)$ independent equations verified by the coefficients $\alpha_{\tau}$. Denote by $F\subset \mathbb{N}^{E(\Gamma)}$ the set of functions $f$ such that $f(e) \in \{0,1\}$ for all $e\in E(\Gamma)$. Lemma~\ref{lemma_meuh} implies that $w_f\in \mathcal{A}_{p,g}\cdot v_0$ for all $ f\in F$. By definition of $v$, we have that $\left(w_f, v \right)_{p,g}^H = 0$ for all $f\in F$. Thus we have
\begin{gather*}\label{eq_cyclicity}
 \sum_{\tau\in [\sigma]}{\left( \prod_{e\in E(\Gamma)}{\lambda_{\tau(e)}^{f(e)}}\right) \alpha_{\tau} ( u_{\tau}, v_0 )_{p,g}^H } =0, \qquad \text{for all} \quad f\in F.
\end{gather*}

Since the complex numbers $ ( u_{\tau}, v_0 )_{p,g}^H$ are nonzero, it is enough to show that the matrix $M := \Big(\prod\limits_{e\in E(\Gamma)}{\lambda_{\tau(e)}^{f(e)}}\Big)_{\substack{\tau\in[\sigma] \\ f\in F}} $ has independent lines to conclude the proof. Recall that we defined Vandermonde matrices $N[e]$ in the proof of Proposition~\ref{prop_cyclicity_high_genus}. Denote by $e_1, \ldots, e_{3g-3}$ the edges of~$\Gamma$. The matrix~$M$ is obtained from the invertible matrix $ N[e_1]\otimes \cdots \otimes N[e_{3g-3}]$ by removing the lines corresponding to non $p$-admissible colorings of $\Gamma$, so has independent lines. This concludes the proof.
\end{proof}

To complete the proof of Proposition~\ref{prop_high_genus}, it remains to show that if $p=4r$, with $r$ an odd prime, verifies Conjecture~\ref{Conjecture} and~$\sigma$ is a coloring of the tetrahedron graph of type~I or~II, then~$u_{\sigma}$ belongs to the cyclic subspace generated by $v_0\in V_p(\Sigma_3)$. We cut the proof in three technical lemmas.

\begin{Lemma}\label{lemma_mama} Suppose $p=4r$ with $r\geq 7$ an odd prime. Let $\sigma$ be a coloration of the tetrahedron graph which is
\begin{enumerate}\itemsep=0pt
\item[$1)$] either of type $I$ and such that either we have $(a\neq \frac{k}{2})$ or $(b\neq \frac{k}{2})$,
\item[$2)$] or of type $II$ and such that either we have $(a\neq b$ and $a\neq k-b)$ or $c\equiv \frac{k}{2} \pmod{4}$,
\end{enumerate}
 then $u_{\sigma}$ belongs to the cyclic space generated by $v_0$.
\end{Lemma}

\begin{proof}The proof is similar in both cases and relies on the following remark: embed a colored tetrahedron graph in $H_3$, choose two opposite edges of the graph and perform two Whitehead moves on these edges. We get this way another embedding of the tetrahedron graph inside $H_3$. While choosing the edges colored by $b$ and its opposite colored by $\frac{k}{2}$ in a type $I$ coloration of the tetrahedron graph, we have
\begin{gather*}
 \Tetahedrongraph{$\frac{k}{2}$}{$\frac{k}{2}$}{$\frac{k}{2}$}{$a$}{$b$}{$c$} = \sum_{i,j} \alpha_{i,j} \invTetahedrongraph{$b$}{$i$}{$\frac{k}{2}$}{$a$}{$j$}{$\frac{k}{2}$},
\end{gather*}
where
\begin{gather*}
\alpha_{i,j}=\left< \Tetahedrongraph{$a$}{$i$}{$\frac{k}{2}$}{$b$}{$\frac{k}{2}$}{$\frac{k}{2}$}\right> \left< \Tetahedrongraph{$\frac{k}{2}$}{$j$}{$\frac{k}{2}$}{$b$}{$c$}{$a$} \right>c_{i,j} \qquad \text{with} \quad
c_{i,j}\neq 0.
\end{gather*}

 If $i=\frac{k}{2}$ then $\alpha_{i,j}=0$. If $i\neq \frac{k}{2}$ then the vector $ \invTetahedrongraph{$b$}{$i$}{$\frac{k}{2}$}{$a$}{$j$}{$\frac{k}{2}$}$ belongs to $Z_{p,3}$, since either $a\neq \frac{k}{2}$ or $b\neq \frac{k}{2}$. Hence this vector belongs to the cyclic space generated by $v_0$ by Lemma~\ref{lemma_ultimate}. This concludes the proof in the first case. The second case is proved similarly.
\end{proof}

\begin{Lemma} Suppose $p=4r$ with $r\geq 7$ an odd prime. Let $\sigma$ be a coloration of the tetrahedron graph of type $I$, such that $a=b= \frac{k}{2}$. Then $u_{\sigma}$ belongs to the cyclic space generated by $v_0$.
\end{Lemma}

\begin{proof}Lemma \ref{fusion_rule} implies that{\samepage
\begin{gather*} v_a:=\Tetahedrongraph{$\frac{k}{2}$}{$a$}{$\frac{k}{2}$}{$\frac{k}{2}$}{$\frac{k}{2}$}{$\frac{k}{2}$} = \left< \thetagraph{$\frac{k}{2}$}{$\frac{k}{2}$}{$\frac{k}{2}$} \right>
\tikz[scale=0.4,baseline=-0.5ex]{
\draw (-2,0) circle (1);
\draw (2,0) circle (1);
\draw (-2,-1) arc (-135: -45: 2.6);
\draw (2,1) [dotted] arc ( 45:135: 2.6);
\draw (0, 1.5) node[above]{$0$};
\draw (0, -1.5) node[below]{$a$};
\draw (-3, 0) node[left]{$\frac{k}{2}$};
\draw (3, 0) node[right]{$\frac{k}{2}$};
} + v',\end{gather*}
 where $v'$ is a vector orthogonal to the first one.}

Proposition \ref{prop_cyclicity_high_genus} implies that the vacuum vector is cyclic in genus $2$ and Lemma~\ref{lemma_reduction} implies that the vector $
\tikz[scale=0.4,baseline=-0.5ex]{
\draw (-2,0) circle (1);
\draw (2,0) circle (1);
\draw (-2,-1) arc (-135: -45: 2.6);
\draw (2,1) [dotted] arc ( 45:135: 2.6);
\draw (0, 1.5) node[above]{$0$};
\draw (0, -1.5) node[below]{$a$};
\draw (-3, 0) node[left]{$\frac{k}{2}$};
\draw (3, 0) node[right]{$\frac{k}{2}$};
}$
belongs to the cyclic space generated by the vacuum vector in genus $3$.

When $a=\frac{k}{2}$, then $W_{[v_a]}$ is one dimensional, so according to Lemma~\ref{lemma1}, either $v_{\frac{k}{2}}$ belongs the cyclic space generated by $v_0$ or it belongs to its orthogonal. Since its scalar product with the vector
$
\tikz[scale=0.4,baseline=-0.5ex]{
\draw (-2,0) circle (1);
\draw (2,0) circle (1);
\draw (-2,-1) arc (-135: -45: 2.6);
\draw (2,1) [dotted] arc ( 45:135: 2.6);
\draw (0, 1.5) node[above]{$0$};
\draw (0, -1.5) node[below]{$a$};
\draw (-3, 0) node[left]{$\frac{k}{2}$};
\draw (3, 0) node[right]{$\frac{k}{2}$};
}$
is a non null $3j$-symbol, the vector $v_a$ is in the cyclic space of~$v_0$.

When $a\neq \frac{k}{2}$, then $W_{[v_a]}$ is two dimensional generated by $v_a$ and $v_{k-a}$. If $v=\alpha_1 v_a + \alpha_2 v_{k-a}$ belongs to the orthogonal of the cyclic space generated by $v_0$, then $v$ is orthogonal to both vectors
$
\tikz[scale=0.4,baseline=-0.5ex]{
\draw (-2,0) circle (1);
\draw (2,0) circle (1);
\draw (-2,-1) arc (-135: -45: 2.6);
\draw (2,1) [dotted] arc ( 45:135: 2.6);
\draw (0, 1.5) node[above]{$0$};
\draw (0, -1.5) node[below]{$a$};
\draw (-3, 0) node[left]{$\frac{k}{2}$};
\draw (3, 0) node[right]{$\frac{k}{2}$};
}$ and
$
\tikz[scale=0.4,baseline=-0.5ex]{
\draw (-2,0) circle (1);
\draw (2,0) circle (1);
\draw (-2,-1) arc (-135: -45: 2.6);
\draw (2,1) [dotted] arc ( 45:135: 2.6);
\draw (0, 1.5) node[above]{$0$};
\draw (0, -1.5) node[below]{$k-a$};
\draw (-3, 0) node[left]{$\frac{k}{2}$};
\draw (3, 0) node[right]{$\frac{k}{2}$};
}$. This implies that $v=0$ so $W_{[v_a]}$ is included in the cyclic space of the vacuum vector.
\end{proof}

\begin{Lemma} Suppose $p=4r$ with $r\geq 7$ an odd prime. Let $\sigma$ be a coloration of the tetrahedron graph of type $II$, such that we have $($either $a= b$ or $a = k-b)$ and $(c\equiv \frac{k}{2}+2 \pmod{4})$. Then~$u_{\sigma}$ belongs to the cyclic space generated by~$v_0$.
\end{Lemma}

\begin{proof}Lemma \ref{fusion_rule} implies that
\begin{equation*} \Tetahedrongraph{$\frac{k}{2}$}{$a$}{$k{-}a$}{$a$}{$c$}{$a$} = \sum_i \sixj{a}{a}{\frac{k}{2}}{k-a}{a}{i}
\tikz[scale=0.4,baseline=-0.5ex]{
\draw (0,1) circle (1);
\draw (0,0) circle (3);
\draw (0,0) --(0,-3);
\draw (0,2)--(0,3);
\draw (0, 2.5) node[left]{$i$};
\draw (-1, 1) node[left]{$a$};
\draw (1, 1) node[right]{$k{-}a$};
\draw (-3, 0) node[left]{$a$};
\draw (0,-2) node[left]{$c$};
\draw (3,0) node[right]{$a$};
}.
\end{equation*}

Let $T\in \emcg$ represents the lift of the Dehn twist around the edge colored by $i$ in the above graph. We have
\begin{equation*} \rho_{p,3}(T)\cdot \Tetahedrongraph{$\frac{k}{2}$}{$a$}{$k{-}a$}{$a$}{$c$}{$a$} = \sum_i \sixj{a}{a}{\frac{k}{2}}{k-a}{a}{i}(-1)^iA^{i(i+2)}
\tikz[scale=0.4,baseline=-0.5ex]{
\draw (0,1) circle (1);
\draw (0,0) circle (3);
\draw (0,0) --(0,-3);
\draw (0,2)--(0,3);
\draw (0, 2.5) node[left]{$i$};
\draw (-1, 1) node[left]{$a$};
\draw (1, 1) node[right]{$k{-}a$};
\draw (-3, 0) node[left]{$a$};
\draw (0,-2) node[left]{$c$};
\draw (3,0) node[right]{$a$};
}.\end{equation*}

Applying the preceding Whitehead move in the opposite direction, we see that $\rho_{p,3}(T)\cdot u_{\sigma}$ belongs to the space generated by the vectors of the form $\Tetahedrongraph{$j$}{$a$}{$k{-}a$}{$a$}{$c$}{$a$}$. Whenever $j\neq \frac{k}{2}$, these generating vectors belong to~$Z_{p,3}$ and thus to the cyclic space generated by the vacuum vector. Denote by $\beta$ the scalar product $\langle u_{\sigma}, \rho_{p,3}(T)u_{\sigma}\rangle$.

If $\beta=0$, then $\rho_{p,3}(T)\cdot u_{\sigma}$ belongs to the cyclic space generated by $v_0$, so does $u_{\sigma}$ since $\rho_{p,3}(T)$ is invertible.

If $\beta\neq 0$, then the operator $A:=\beta \cdot \mathds{1} + \sixj{a}{a}{\frac{k}{2}}{k-a}{a}{\frac{k}{2}} \rho_{p,3}(T) \in \mathcal{A}_{p,3}$ is invertible since~$\rho_{p,3}(T)$ has finite order. Since $A\cdot u_{\sigma}$ belongs to the cyclic space generated by the vacuum vector, so does~$u_{\sigma}$.
\end{proof}

\appendix

\section[Appendix: Vanishing $6j$-symbols]{Appendix: Vanishing $\boldsymbol{6j}$-symbols}\label{appendixA}
Assume that $p=2(k+2)$ with $k\geq 2$ an even integer. The goal of this appendix is to prove Proposition~\ref{prop_null_6j}.

\begin{Lemma}\label{lemma_A}Let $a$, $b$ be two integers such that $(a,k-a,b)$ is $p$-admissible. Define
\begin{gather*} F(a,b):=\frac{\left< \Tetahedrongraph{a}{k}{$k{-}a$}{$k{-}a$}{b}{a} \right> \langle a \rangle }{\left<\thetagraph{$k{-}a$}{a}{b}\right>}.\end{gather*}
Then we have $F(a,b)=(-1)^{\frac{b+k}{2}+a}$.
\end{Lemma}

\begin{proof}A straightforward computation using the formulas of \cite{MV} gives $F(a,b)= \frac{f(a)}{g(b)}$, where $f(a)=(-1)^a [a+1]![k-a]!$ and $g(b)=(-1)^{\frac{k+b}{2}}\left[\frac{k-b}{2}\right]!\left[\frac{k+b}{2}+1\right]! $. Note that
\begin{gather*}
\frac{f(a+1)}{f(a)}=-\frac{[a+2]}{[k-a]}=-1\qquad \text{and} \qquad \frac{g(b+2)}{g(b)}=-\frac{\left[\frac{k+b}{2}+2\right]}{\left[\frac{k-b}{2}\right]}=-1.\end{gather*}
We conclude using the fact that $F\big(\frac{k}{2},2\big)=-1$.
\end{proof}

\begin{Lemma}\label{lemma_B}If $(a,b,c)$ is a $p$-admissible triple, then we have
\begin{gather*} \langle k-a\rangle \langle k-b\rangle \frac{\left< \Tetahedrongraph{a}{k}{b}{$k{-}a$}{c}{$k{-}b$}\right>}{\left< \thetagraph{$k{-}a$~~}{$k{-}b$}{c} \right>} \cdot \frac{\left< \Tetahedrongraph{$k{-}a$}{k}{$k{-}b$}{a}{c}{b}\right>}{\left< \thetagraph{a~~}{$k{-}b$}{c} \right>}=1.
\end{gather*}
\end{Lemma}

\begin{proof}We use the fact that adding a trivial framed knot colored by $k$ does not change the class of a vector. We work in the space associated to the sphere with three punctures colored by $a$, $b$ and $c$ and compute
\begin{gather*}
\tikz[scale=0.3,baseline=-0.5ex]{
\draw (0,0)--(-3, 3);
\draw (0,0)--(3, 3);
\draw (0,-4)--(0,0);
\draw (-2,3) node[above]{$a$};
\draw (2,3) node[above]{$b$};
\draw (0,-3) node[right]{$c$};
}
=
\tikz[scale=0.3,baseline=-0.5ex]{
\draw (0,2) circle (1);
\draw (0,0)--(-3, 3);
\draw (0,0)--(3, 3);
\draw (0,-4)--(0,0);
\draw (-2,3) node[above]{$a$};
\draw (2,3) node[above]{$b$};
\draw (0,-3) node[right]{$c$};
\draw (0,3) node[above]{$k$};
} = \langle k-a\rangle \langle k-b\rangle
\tikz[scale=0.3,baseline=-0.5ex]{
\draw (0,0)--(-3, 3);
\draw (0,0)--(3, 3);
\draw (0,-4)--(0,0);
\draw (-1.5, 1.5)--(1.5,1.5);
\draw (-2.7,2.7)--(2.7,2.7);
\draw (-1,0.5) node[left]{$a$};
\draw (1, 0.5) node[right]{$b$};
\draw (0, 3) node[above]{$k$};
\draw (0, 1.5) node[above]{$k$};
\draw (-2,2) node[left]{$k-a$};
\draw (2,2) node[right]{$k-b$};
\draw (-3,4) node[above]{$a$};
\draw (3,4) node[above]{$b$};
\draw (0,-3) node[right]{$c$};
}
\\
 \hphantom{\tikz[scale=0.3,baseline=-0.5ex]{
\draw (0,0)--(-3, 3);
\draw (0,0)--(3, 3);
\draw (0,-4)--(0,0);
\draw (-2,3) node[above]{$a$};
\draw (2,3) node[above]{$b$};
\draw (0,-3) node[right]{$c$};
}}{}
 =
\langle k-a\rangle \langle k-b\rangle \frac{\left< \Tetahedrongraph{$a$}{$k$}{$b$}{$k{-}a$}{$c$}{$k{-}b$}\right>}{\left< \thetagraph{$k{-}a~$}{$k{-}b$}{$c$} \right>} \cdot \frac{\left< \Tetahedrongraph{$k{-}a$}{$k$}{$k{-}b$}{$a$}{$c$}{$b$}\right>}{\left< \thetagraph{$a~$}{$k{-}b$}{$c$} \right>}
\tikz[scale=0.3,baseline=-0.5ex]{
\draw (0,0)--(-3, 3);
\draw (0,0)--(3, 3);
\draw (0,-4)--(0,0);
\draw (-2,3) node[above]{$a$};
\draw (2,3) node[above]{$b$};
\draw (0,-3) node[right]{$c$};
}.
\end{gather*}
We conclude the proof by identifying both vectors.
\end{proof}

\begin{proof}[Proof of Proposition \ref{prop_null_6j}] We use the fact that the Kauffman bracket of a link in $S^3$ does not change if we add a trivial framed knot colored by~$k$. First, for a $6j$-symbol of type~$I$, we use Lemmas~\ref{fusion_rule} and~\ref{lemma_A} to compute
\begin{gather*}
\tikz[scale=0.4,baseline=-0.5ex]{
\draw (0,0) circle (2);
\draw (0,1) circle (0.6);
\draw (0,0)--(0,-2);
\draw (-1.35,1.35)--(0,0);
\draw (1.35,1.35)--(0,0);
\draw (0,2) node[above]{$\frac{k}{2}$};
\draw (-1,1) node[below]{$\frac{k}{2}$};
\draw (1,1) node[below]{$\frac{k}{2}$};
\draw (-1.5,-1.3) node[left]{$a$};
\draw (0,-1.3) node[left]{$b$};
\draw (1.5,-1.3) node[right]{$c$};
\draw (0,0.6) node[above]{$k$};
} =
 \left< \tfrac{k}{2}\right> ^3
\tikz[scale=0.4,baseline=-0.5ex]{
\draw (-1,0)--(1,0);
\draw (-1,0)--(0,-1);
\draw (1,0)--(0,-1);
\draw (-1,0)--(-3,1);
\draw (-3,1)--(-3,3);
\draw (-3,1)--(-4,2);
\draw (-4,2)--(-3,3);
\draw (1,0)--(3,1);
\draw (3,1)--(3,3);
\draw (3,1)--(4,2);
\draw (4,2)--(3,3);
\draw (0,-1)--(0,-2);
\draw (-4,2) arc (-180:0:4);
\draw (3,3) arc (20:160:3.2);
\draw (0, -1.5) node[left]{$b$};
\draw (-3,-1) node[left]{$a$};
\draw (3,-1) node[right]{$c$};
\draw (-0.5,-0.5) node[left]{$\frac{k}{2}$};
\draw (0.5,-0.5) node[right]{$\frac{k}{2}$};
\draw (0,0) node[above]{$k$};
\draw (-2,0.5) node[above]{$\frac{k}{2}$};
\draw (2,0.5) node[above]{$\frac{k}{2}$};
\draw (-3.5, 1.5) node[below]{$\frac{k}{2}$};
\draw (3.5, 1.5) node[below]{$\frac{k}{2}$};
\draw (-3.2,2) node[right]{$k$};
\draw (3.2,2) node[left]{$k$};
\draw (3.5,2.5)node[above]{$\frac{k}{2}$};
\draw (-3.5,2.5)node[above]{$\frac{k}{2}$};
\draw (0,5) node[above]{$\frac{k}{2}$};
} \\
\hphantom{\tikz[scale=0.4,baseline=-0.5ex]{
\draw (0,0) circle (2);
\draw (0,1) circle (0.6);
\draw (0,0)--(0,-2);
\draw (-1.35,1.35)--(0,0);
\draw (1.35,1.35)--(0,0);
\draw (0,2) node[above]{$\frac{k}{2}$};
\draw (-1,1) node[below]{$\frac{k}{2}$};
\draw (1,1) node[below]{$\frac{k}{2}$};
\draw (-1.5,-1.3) node[left]{$a$};
\draw (0,-1.3) node[left]{$b$};
\draw (1.5,-1.3) node[right]{$c$};
\draw (0,0.6) node[above]{$k$};
}}{}
=
F \big(\tfrac{k}{2}, a \big)F\big(\tfrac{k}{2}, b \big)F\big(\tfrac{k}{2}, c \big) \Tetahedrongraph{$\frac{k}{2}$}{$\frac{k}{2}$}{$\frac{k}{2}$}{$a$}{$b$}{$c$}
= -\Tetahedrongraph{$\frac{k}{2}$}{$\frac{k}{2}$}{$\frac{k}{2}$}{$a$}{$b$}{$c$}.
\end{gather*}

Thus $\left< \Tetahedrongraph{$\frac{k}{2}$}{$\frac{k}{2}$}{$\frac{k}{2}$}{$a$}{$b$}{$c$} \right> =0$.

For a $6j$-symbol of type $II$, a similar computation using Lemmas \ref{lemma_A} and \ref{lemma_B} gives
\begin{gather*}
\tikz[scale=0.4,baseline=-0.5ex]{
\draw (0,0) circle (2);
\draw (0,1) circle (0.6);
\draw (0,0)--(0,-2);
\draw (-1.35,1.35)--(0,0);
\draw (1.35,1.35)--(0,0);
\draw (0,2) node[above]{$\frac{k}{2}$};
\draw (-1,1) node[below]{$a$};
\draw (1.2,1) node[below]{$k{-}a$};
\draw (-1.5,-1.3) node[left]{$b$};
\draw (0,-1.3) node[left]{$c$};
\draw (1.5,-1.3) node[right]{$b$};
\draw (0,0.6) node[above]{$k$};
} = \left< \tfrac{k}{2}\right> \langle a \rangle \langle k-a \rangle
\tikz[scale=0.4,baseline=-0.5ex]{
\draw (-1,0)--(1,0);
\draw (-1,0)--(0,-1);
\draw (1,0)--(0,-1);
\draw (-1,0)--(-3,1);
\draw (-3,1)--(-3,3);
\draw (-3,1)--(-4,2);
\draw (-4,2)--(-3,3);
\draw (1,0)--(3,1);
\draw (3,1)--(3,3);
\draw (3,1)--(4,2);
\draw (4,2)--(3,3);
\draw (0,-1)--(0,-2);
\draw (-4,2) arc (-180:0:4);
\draw (3,3) arc (20:160:3.2);
\draw (0, -1.5) node[left]{$c$};
\draw (-3,-1) node[left]{$b$};
\draw (3,-1) node[right]{$b$};
\draw (-0.5,-0.5) node[left]{$a$};
\draw (0.5,-0.5) node[right]{$k{-}a$};
\draw (0,0) node[above]{$k$};
\draw (-2,0.5) node[above]{$k{-}a$};
\draw (2,0.5) node[above]{$a$};
\draw (-3.5, 1.5) node[below]{$a$};
\draw (3.6, 1.5) node[below]{$k{-}a$};
\draw (-3,2) node[right]{$k$};
\draw (3,2) node[left]{$k$};
\draw (3.5,2.5)node[above]{$\frac{k}{2}$};
\draw (-3.5,2.5)node[above]{$\frac{k}{2}$};
\draw (0,5) node[above]{$\frac{k}{2}$};
} \\
=
F (a, c ) \cdot \left(
\left<\tfrac{k}{2} \right>\langle k-a\rangle \frac{\left< \Tetahedrongraph{$\frac{k}{2}$}{$k$}{$a$}{$\frac{k}{2}$}{$b$}{$k{-}a$}\right>}{\left< \thetagraph{$k{-}a$}{$\frac{k}{2}$}{$b$} \right>} \cdot \frac{\left< \Tetahedrongraph{$\frac{k}{2}$}{$k$}{$k{-}a$}{$\frac{k}{2}$}{$b$}{$a$}\right>}{\left< \thetagraph{$k{-}a$}{$\frac{k}{2}$}{$b$} \right>}
\right) \Tetahedrongraph{$\frac{k}{2}$}{$a$}{$k{-}a$}{$b$}{$c$}{$b$}
\\ = -\Tetahedrongraph{$\frac{k}{2}$}{$a$}{$k{-}a$}{$b$}{$c$}{$b$}.
\end{gather*}
Thus $\left< \Tetahedrongraph{$\frac{k}{2}$}{$a$}{$k{-}a$}{$b$}{$c$}{$b$} \right> =0$.
\end{proof}

\subsection*{Acknowledgements}

The author is thankful to L.~Funar and F.~Costantino for useful discussions
and to the anonymous referees for valuable corrections that improved the
clarity of the paper. He acknowledges support from the grant ANR $2011$ BS $01 020 01$ ModGroup, CAPES and the GEAR Network.

\pdfbookmark[1]{References}{ref}
\LastPageEnding


\begin{thebibliography}{99}
\footnotesize\itemsep=0pt

\bibitem{AF}
Andersen J.E., Fjelstad J., Reducibility of quantum representations of mapping
 class groups, \href{https://doi.org/10.1007/s11005-009-0367-7}{\textit{Lett. Math. Phys.}} \textbf{91} (2010), 215--239,
 \href{https://arxiv.org/abs/0806.2539}{arXiv:0806.2539}.

\bibitem{BHMV2}
Blanchet C., Habegger N., Masbaum G., Vogel P., Topological quantum field
 theories derived from the {K}auffman bracket, \href{https://doi.org/10.1016/0040-9383(94)00051-4}{\textit{Topology}} \textbf{34}
 (1995), 883--927.

\bibitem{FK}
Freedman M., Krushkal V., On the asymptotics of quantum {${\rm SU}(2)$}
 representations of mapping class groups, \href{https://doi.org/10.1515/FORUM.2006.017}{\textit{Forum Math.}} \textbf{18}
 (2006), 293--304, \href{https://arxiv.org/abs/math.QA/0409503}{arXiv:math.QA/0409503}.

\bibitem{GM13}
Gilmer P.M., Masbaum G., Maslov index, lagrangians, mapping class groups and
 {TQFT}, \href{https://doi.org/10.1515/form.2011.143}{\textit{Forum Math.}} \textbf{25} (2013), 1067--1106,
 \href{https://arxiv.org/abs/0912.4706}{arXiv:0912.4706}.

\bibitem{Jones_proj}
Jones V.F.R., Index for subfactors, \href{https://doi.org/10.1007/BF01389127}{\textit{Invent. Math.}} \textbf{72} (1983),
 1--25.

\bibitem{KoberdaSantharoubane17}
Koberda T., Santharoubane R., Irreducibility of quantum representations of
 mapping class groups with boundary, \href{https://doi.org/10.4171/QT/116}{\textit{Quantum Topol.}} \textbf{9}
 (2018), 633--641, \href{https://arxiv.org/abs/1701.08901}{arXiv:1701.08901}.

\bibitem{Koju1}
Korinman J., Irreducible factors of {W}eil representations and {TQFT},
 \textit{Math. Rep.}, {t}o appear, \href{https://arxiv.org/abs/1310.0390}{arXiv:1310.0390}.

\bibitem{Li2}
Lickorish W.B.R., Invariants for 3-manifolds from the combinatorics of the
 {J}ones polynomial, \href{https://doi.org/10.2140/pjm.1991.149.337}{\textit{Pacific~J. Math.}} \textbf{149} (1991), 337--347.

\bibitem{MR}
Masbaum G., Roberts J.D., On central extensions of mapping class groups,
 \href{https://doi.org/10.1007/BF01444490}{\textit{Math. Ann.}} \textbf{302} (1995), 131--150, \href{https://arxiv.org/abs/math.QA/9909128}{arXiv:math.QA/9909128}.

\bibitem{MV}
Masbaum G., Vogel P., {$3$}-valent graphs and the {K}auffman bracket,
 \href{https://doi.org/10.2140/pjm.1994.164.361}{\textit{Pacific~J. Math.}} \textbf{164} (1994), 361--381.

\bibitem{RT}
Reshetikhin N., Turaev V.G., Invariants of {$3$}-manifolds via link polynomials
 and quantum groups, \href{https://doi.org/10.1007/BF01239527}{\textit{Invent. Math.}} \textbf{103} (1991), 547--597.

\bibitem{Ro}
Roberts J., Irreducibility of some quantum representations of mapping class
 groups, \href{https://doi.org/10.1142/S021821650100113X}{\textit{J.~Knot Theory Ramifications}} \textbf{10} (2001), 763--767,
 \href{https://arxiv.org/abs/math.QA/9909128}{arXiv:math.QA/9909128}.

\bibitem{Wenzl}
Wenzl H., On sequences of projections, \textit{C.~R.~Math. Rep. Acad. Sci.
 Canada} \textbf{9} (1987), 5--9.

\bibitem{Wi2}
Witten E., Quantum field theory and the {J}ones polynomial, \href{http://dx.doi.org/10.1007/BF01217730}{\textit{Comm. Math.
 Phys.}} \textbf{121} (1989), 351--399.

\end{thebibliography}
\end{document}